\documentclass[10pt]{article}

\usepackage{amsfonts} 
\usepackage{amsmath,hhline,latexsym}
\usepackage{amssymb}

\usepackage[latin1]{inputenc}

\usepackage{hyperref}

\usepackage{color}
\usepackage{graphicx}

\newtheorem{hyp}{Hypothesis}

\newtheorem{theorem}{\sc Theorem}[section]
\newtheorem{lemma}{\sc Lemma}[section]
\newtheorem{proposition}{\sc Proposition}[section]

\newtheorem{remark}{Remark}

\newcommand{\dis}{\displaystyle}
\newcommand{\eps}{\varepsilon}
\newcommand{\om}{\omega}

\newcommand{\ph}{\varphi}

\newcommand{\Fin}{\hfill$\Box$}

\newcommand{\N}{\mbox{$I \kern -4pt N$}}
\newcommand{\Q}{\mbox{$Q \kern -8pt I$}}
\newcommand{\R}{\mbox{$I \kern -4pt R$}}
\newcommand{\C}{\mbox{$C \kern -8pt I$}}

\newcommand{\jnt}{\dis\int}
\newcommand{\jjntQT}{\jnt\!\!\!\!\jnt_{Q_{T}}}
\newcommand{\jjntqT}{\jnt\!\!\!\!\jnt_{q_{T}}}

\providecommand{\tabularnewline}{\\}

\def\R{{\bf R}}
\title{\textbf{Inverse problems for linear hyperbolic equations using mixed formulations}}

\author{ 
	\textsc{Nicolae C\^indea}\thanks{Laboratoire de Math\'ematiques, Universit\'e Blaise Pascal (Clermont-Ferrand 2), UMR CNRS 6620, 
	Campus de C\'ezeaux, 63177, Aubi\`ere, France. E-mails: {\tt nicolae.cindea@math.univ-bpclermont.fr.}}\quad 
	\and
	\textsc{Arnaud M\"unch}\thanks{Laboratoire de Math\'ematiques, Universit\'e Blaise Pascal (Clermont-Ferrand 2), UMR CNRS 6620, 
	Campus de C\'ezeaux, 63177, Aubi\`ere, France.
	E-mail: {\tt arnaud.munch@math.univ-bpclermont.fr}.}}
	
\begin{document}

\maketitle

\begin{abstract}
We introduce in this document a direct method allowing to solve numerically inverse type problems for linear hyperbolic equations. We first consider the reconstruction of the full solution of the wave equation posed in $\Omega\times (0,T)$ - $\Omega$ a bounded subset of $\mathbb{R}^N$ - from a partial distributed observation. We employ a least-squares technique and minimize the $L^2$-norm of the distance from the observation to any solution. Taking the hyperbolic equation as the main constraint of the problem, the optimality conditions are reduced to a mixed formulation involving both the state to reconstruct and a Lagrange multiplier. Under usual geometric optic conditions, we show the well-posedness of this mixed formulation (in particular the inf-sup condition) and then introduce a numerical approximation based on space-time finite elements discretization. We prove the strong convergence of the approximation and then discussed several examples for $N=1$ and $N=2$. The problem of the reconstruction of both the state and the source term is also addressed.


\end{abstract}

\textsc{keywords :}Linear wave equation, Inverse problem, Finite elements methods; Mixed formulation

\textsc{AMS number:} 35L10, 65M12,  93B40.

\section{Introduction - Inverse problems for hyperbolic equations}
\label{sec:intro}

Let $\Omega$ be a bounded domain of $\mathbb{R}^N$ ($N\geq 1$) whose boundary $\partial\Omega$ is Lipschitz and let $T>0$. We note $Q_T:=\Omega\times (0,T)$ and $\Sigma_T:=\partial\Omega\times (0,T)$. We are concerned in this work with inverse type problems for linear hyperbolic equation of the following type 
\begin{equation}
 \label{eq:wave}
 \left\{
   \begin{array}{ll}
   y_{tt} - \nabla\cdot (c(x) \nabla y) + d(x, t) y= f, & \qquad (x, t) \in Q_T \\
   y = 0, & \qquad (x,t)\in \Gamma_T \\
   (y(\cdot, 0), y_t(\cdot,0)) = (y_0, y_1), & \qquad x \in \Omega.
   \end{array} 
 \right.
\end{equation}
We assume that $c\in C^1(\overline{\Omega}, \mathbb{R})$ with $c(x)\geq c_0>0$ in~$\overline{\Omega}$, $d \in L^\infty(Q_T)$, $(y_0,y_1) \in \boldsymbol{H}:=L^2(\Omega)\times H^{-1}(\Omega)$ and $f\in X:=L^2(0,T; H^{-1}(\Omega))$.

   For any $(y_0,y_1)\in \boldsymbol{H}$ and any $f\in X$, there exists exactly one solution $y$ to \eqref{eq:wave}, with $y \in C^0([0, T]; L^2(\Omega)) \cap C^1([0, T]; H^{-1}(\Omega))$ (see~\cite{JLL88}). 
   
   In the sequel, for simplicity, we shall use the following notation:
\begin{equation}
\label{eq:L}
L\, y:=y_{tt}-\nabla\cdot (c(x) \nabla y) + d(x,t)y.
\end{equation}
and $X^{\prime}:=L^2(0,T; H^1_0(\Omega))$.

Let now $\omega$ be any non empty open subset of $\Omega$ and let $q_T:=\omega \times (0,T)\subset Q_T$.
A typical inverse problem for (\ref{eq:wave}) is the following one :  from an \textit{observation} or \textit{measurement} $y_{obs}$ in $L^2(q_T)$ on the sub-domain $q_T$, we want to recover a solution $y$ of the boundary value problem (\ref{eq:wave}) which coincides with the observation on $q_T$. 

Introducing the operator $P:L^2(Q_T)\to X\times L^2(q_T)$ defined by $P\,y:=(Ly,y_{\vert q_T})$, the problem is reformulated as : 
\begin{equation}
\label{IP}\tag{$IP$}
\text{\it find } y\in L^2(Q_T) \text{ \it solution of } P\,y=(f,y_{obs}). 
\end{equation}

From the unique continuation property for (\ref{eq:wave}), if the set $q_T$ satisfies some geometric conditions and if $y_{obs}$ is a restriction to $q_T$ of a solution of (\ref{eq:wave}), then the problem is well-posed in the sense that the state $y$ corresponding to the pair $(y_{obs},f)$ is unique.



In view of the unavoidable uncertainties on the data $y_{obs}$ (coming from measurements, numerical approximations, etc), the problem needs to be relaxed. In this respect, the most natural (and widely used in practice) approach consists to introduce the following extremal problem (of least-squares type)
\begin{equation}
\label{extremal_problem} \tag{\it LS}
\quad
\left\{
\begin{aligned}
&\textrm{minimize over } \boldsymbol{H} \quad J(y_0,y_1):=\frac{1}{2}  \Vert y-y_{obs} \Vert^2_{L^2(q_T)}\\
& \textrm{where} \quad y \quad \textrm{solves} \quad (\ref{eq:wave}),
\end{aligned}
\right.
\end{equation}
since $y$ is uniquely and fully determined from $f$ and the data $(y_0,y_1)$. Here the constraint $y-y_{obs}=0$ in $L^2(q_T)$ is relaxed; however, if $y_{obs}$ is a restriction to $q_T$ of a solution of (\ref{eq:wave}), then problems \eqref{extremal_problem} and \eqref{IP} obviously coincide. A minimizing sequence for $J$ in $\boldsymbol{H}$ is easily defined in term of the solution of an auxiliary adjoint problem. Apart from a possible low decrease of the sequence near extrema, the main drawback, when one wants to prove the convergence of a discrete approximation is that, it is in general not possible to minimize over a discrete subspace of $\{y; Ly-f=0\}$ subject to the equality (in $X$) $Ly-f=0$. Therefore, the minimization procedure first requires the discretization of the functional $J$ and of the system (\ref{eq:wave}); this raised the issue of uniform coercivity property (typically here some uniform discrete observability inequality for the adjoint solution) of the discrete functional with respect to the approximation parameter. As far as we know, this delicate issue has received answers only for specific and somehow academic situations (uniform Cartesian approximation of $\Omega$, constant coefficients in \eqref{eq:wave}). We refer to \cite{NC-AM-mixedwave,Glo08,komornikloreti,munch05} and the references therein.

More recently, a different method to solve inverse type problems like \eqref{IP} has emerged and use so called Luenberger type observers:  this consists in defining, from the observation on $q_T$, an auxiliary boundary value problem whose solution possesses the same asymptotic behavior in time than the solution of (\ref{eq:wave}): the use of the reversibility of the hyperbolic equation then allows to reconstruct the initial data $(y_0,y_1)$. We refer to \cite{cindea_moireau,ramdani2010} and the references therein. But, for the same reasons, on a numerically point of view, these method require to prove uniform discrete observability properties. 

In a series of works, Klibanov and co-workers use different approaches to solve inverse problems (we refer to \cite{Klibanov-book} and the references therein): they advocate in particular the quasi-reversibility method which reads as follows : for any $\eps>0$, find $y_{\eps}\in \mathcal{A}$ the solution of 
\begin{equation}\label{FV}\tag{$QR$}
 \langle Py_\eps,P\overline{y} \rangle_{X\times L^2(q_T)} + \eps \langle y_{\eps},\overline{y} \rangle_{\mathcal{A}} = \left\langle(f,y_{obs}), P\overline{y}\right\rangle_{X^{\prime}\times L^2(q_T),X\times L^2(q_T)},
\end{equation}
for all $\overline{y}\in \mathcal{A}$, 
where $\mathcal{A}$ denotes a Hilbert space subset of $L^2(Q_T)$ so that $Py\in X\times L^2(q_T)$ for all $y\in \mathcal{A}$ and $\eps>0$ a Tikhonov like parameter which ensures the well-posedness. We refer for instance to \cite{clason_klibanov} where the lateral Cauchy problem for the wave equation with non constant diffusion is addressed within this method. Remark that (\ref{FV}) can be viewed as a least-squares problem since the solution $y_{\eps}$ minimizes over $\mathcal{A}$ the functional $y\to \Vert P y - (f,y_{obs})\Vert^2_{X\times L^2(q_T)}+ \eps \Vert y\Vert^2_{\mathcal{A}}$. Eventually, if $y_{obs}$ is a restriction to $q_T$ of a solution of (\ref{eq:wave}), 
the corresponding $y_{\eps}$ converges in $L^2(Q_T)$ toward to the solution of \eqref{IP} as $\eps\to 0$. There, unlike in Problem \eqref{extremal_problem}, the unknown is the state variable $y$ itself (as it is natural for elliptic equations) so that any standard numerical methods based on a conformal approximation of the space $\mathcal{A}$ together with appropriate observability inequalities allow to obtain a convergent approximation of the solution. In particular, there is no need to prove discrete observability inequalities. We refer to the book \cite{BeilinaKlibanov14}. We also mention \cite{bourgeois2010,bourgeoisbis} where a similar technique has been used recently to solve the inverse obstacle problem associated to the Laplace equation, which consists in finding an interior obstacle from boundary Cauchy data. 

In the spirit of the works \cite{Klibanov-book,bourgeois2010,clason_klibanov}, we explore the direct resolution of the optimality conditions associated to the extremal problem \eqref{extremal_problem}, without Tikhonov parameter while keeping $y$ as the unknown of the problem. This strategy, which avoids any iterative process, has been successfully applied in the closed context of the exact controllability of (\ref{eq:wave}) in \cite{NC-AM-mixedwave} and \cite{CC-NC-AM,NC-EFC-AM}. The idea is to take into account the state constraint $Ly-f=0$ with a Lagrange multiplier. This allows to derive explicitly the optimality systems associated  to \eqref{extremal_problem} in term of an elliptic mixed formulation and therefore reformulate the original problem. Well-posedness of such new formulation is related to an observability inequality for the homogeneous solution of the hyperbolic equation. 

The outline of this paper is as follow. In Section \ref{recovering_y}, we consider the least-squares problem \eqref{P} and reconstruct the solution of the wave equation from a partial observation localized on a subset $q_T$ of $Q_T$. For that, in Section \ref{sec2_direct}, we associate to \eqref{P} the equivalent mixed formulation (\ref{eq:mf}) which relies on the optimality conditions of the problem. Assuming that $q_T$ satisfies the classical geometric optic condition (Hypothesis 1, see (\ref{iobs})), we then show the well-posedness of this mixed formulation, in particular, we check the Babuska-Brezzi inf-sup condition (see Theorem \ref{th:mf}). Interestingly, in Section \ref{sec2_dual}, we also derive a equivalent dual extremal problem, which reduces the determination of the state $y$ to the minimization of an elliptic functional with respect to the Lagrange multiplier. In Section \ref{recovering_y_f}, we apply the same procedure to recover from a partial observation both the state and the source term. Section \ref{sec_numer} is devoted to the numerical approximation, through a conformal space-time finite element discretization. The strong convergence of the approximation $(y_h,f_h)$ is shown as the discretization parameter $h$ tends to zero. In particular, we discuss the discrete  inf-sup property of the mixed formulation. We present numerical experiments in Section \ref{sec_experiment} for $\Omega=(0,1)$ and $\Omega \subset\mathbb{R}^2$, in agreement with the theoretical part. We consider in particular time dependent observation zones. Section \ref{sec_conclusion} concludes with some perspectives. 


\section[Recovering the solution from a partial observation]{Recovering the solution from a partial observation: a mixed re-formulation of the problem}\label{recovering_y}

In this section, assuming that the initial $(y_0,y_1)\in \boldsymbol{H}$ are unknown, we address the inverse problem \eqref{IP}. Without loss of generality, in view of the linearity of the system (\ref{eq:wave}), we assume that the source term $f\equiv 0$.  

\par\noindent
We consider the non empty vectorial space $Z$ defined by 
\begin{equation}
Z:=\{y: y\in C([0,T], L^2(\Omega))\cap C^1([0,T], H^{-1}(\Omega)), Ly\in X\}.  \label{refZ}
\end{equation}
and then introduce the following hypothesis :

\begin{hyp}
There exists a constant $C_{obs}=C(\omega,T,\Vert c\Vert_{C^1(\overline{\Omega})},\Vert d\Vert_{L^{\infty}(\Omega)})$ such that the following estimate holds : 
\begin{equation}
\Vert y(\cdot,0),y_t(\cdot,0)\Vert^2_{\boldsymbol{H}} \leq C_{obs} \biggl( \Vert y\Vert^2_{L^2(q_T)}+ \Vert Ly\Vert^2_X  \biggr), \quad \forall y\in Z.  \label{iobs}\tag{$\mathcal{H}$}
\end{equation}
\end{hyp}

Condition \eqref{iobs} is a generalized observability inequality for the solution of the hyperbolic equation: for constant coefficients, this estimate is known to hold if the triplet ($\omega,T,\Omega$) satisfies a geometric optic condition. We refer to \cite{BLR}. In particular, $T$ should be large enough. Upon the same condition, (\ref{iobs}) also holds in the non-cylindrical situation where the domain $\omega$ varies with respect to the time variable: we refer to \cite{CC-NC-AM} for the one dimensional case. For non constant velocity $c$ and potential $d$, we refer to \cite{NC-EFC-AM} and the references therein.

Then, within this hypothesis, for any $\eta>0$, we define on $Z$ the bilinear form 
\begin{equation}
\label{eq:pseta}
\begin{aligned}
\langle y,\overline{y}\rangle_Z:= & \jjntqT y\,\overline{y}\,dxdt + \eta\int_0^T \langle Ly,\, L\overline{y} \rangle_{H^{-1}(\Omega)}\, dt  \quad \forall y,\overline{y}\in Z.
\end{aligned}
\end{equation}
In view of (\ref{iobs}), this bilinear form defines a scalar product over $Z$. Moreover, endowed to this scalar product, we easily obtain that $Z$ is a Hilbert space (see \cite{CC-NC-AM}, Corollary 2.4). We note the corresponding norm by $\Vert y\Vert_{Z}:=\sqrt{\langle y,y \rangle_Z}$.

Then, we consider the following extremal problem : 
\begin{equation}
\label{P}
\tag{$\mathcal{P}$}
\left\{
\begin{aligned}
& \inf  J(y):=  \frac{1}{2}\Vert y-y_{obs}\Vert^2_{L^2(q_T)}, \\
& \textrm{subject to}\quad y\in W
\end{aligned}
\right.
\end{equation}
where $W$ is the closed subspace of $Z$ defined by 
$$
W:=\{y\in Z; \,  Ly=0 \,\, \textrm{in}\,\, X\}
$$
and endowed with the norm of $Z$. 

The extremal problem \eqref{P} is well posed : the functional $J$ is continuous over $W$, is strictly convex and is such that $J(y)\to +\infty$ as $\Vert y\Vert_W\to \infty$. Note also that the solution of \eqref{P} in $W$ does not depend on $\eta$. 

Remind that from the definition of $Z$, $Ly$ belongs to $X$. Similarly, the uniqueness of the solution is lost if the hypothesis \eqref{iobs} is not fulfilled, for instance if $T$ is not large enough.  Eventually, from (\ref{iobs}), the solution $y$ in $Z$ of \eqref{P} satisfies $(y(\cdot,0),y_t(\cdot,0))\in \boldsymbol{H}$, so that problem \eqref{P} is equivalent to the minimization of $J$ with respect to $(y_0,y_1)\in \boldsymbol{H}$ as in problem \eqref{IP}, Section 1. 

We also recall that for any $z\in Z$ there exists a positive constant $C_{\Omega,T}$ such that 
\begin{equation} 
\Vert z\Vert^2_{L^2(Q_T)} \leq C_{\Omega,T} \biggl(  \Vert z(\cdot,0),z_t(\cdot,0)\Vert^2_{\boldsymbol{H}} + \Vert Lz\Vert^2_X\biggr).  \label{continuity_constant}
\end{equation}
This equality and (\ref{iobs}) imply that
\begin{equation}
\Vert z\Vert^2_{L^2(Q_T)} \leq C_{\Omega,T} \biggl(C_{obs} \Vert z\Vert^2_{L^2(q_T)} + (1+C_{obs}) \Vert Lz\Vert^2_X\biggr), \quad \forall z\in Z.    \label{estimate_qTQT}
\end{equation}

\subsection{Direct approach} \label{sec2_direct}

In order to solve \eqref{P}, we have to deal with the constraint equality which appears in the space $W$. Proceeding as in \cite{NC-AM-mixedwave}, we introduce a Lagrangian multiplier $\lambda\in X^{\prime}$ and the following mixed formulation: find $(y, \lambda)\in Z\times X^{\prime}$ solution of 
\begin{equation} \label{eq:mf}
\left\{
\begin{array}{rcll}
\noalign{\smallskip} a(y, \overline{y}) + b(\overline{y}, \lambda) & = & l(\overline{y}), & \qquad \forall \overline{y} \in Z \\
\noalign{\smallskip} b(y, \overline{\lambda}) & = & 0, & \qquad \forall \overline{\lambda} \in X^{\prime},
\end{array}
\right.
\end{equation}
where
\begin{align}
\label{eq:a} & a : Z \times Z \to \mathbb{R},  \quad a(y,\overline{y}) := \jjntqT y\,\overline{y}\,dxdt,
\\
\label{eq:b} & b: Z \times X^{\prime}  \to \mathbb{R},  \quad b(y,\lambda) := \int_0^T \langle \lambda,\, Ly \rangle_{H^1_0(\Omega),H^{-1}(\Omega)} dt,\\
\label{eq:l} & l: Z \to \mathbb{R},  \quad l(y) := \jjntqT y_{obs} \,y\, dxdt. 
\end{align}
System (\ref{eq:mf}) is nothing else than the optimality system corresponding to the extremal problem \eqref{P}. Precisely, the following result holds : 
\begin{theorem}\label{th:mf} Under the hypothesis \eqref{iobs}, 
\begin{enumerate}
\item The mixed formulation (\ref{eq:mf}) is well-posed.
\item The unique solution $(y, \lambda) \in Z\times X^{\prime}$ to \eqref{eq:mf} is the unique saddle-point of the Lagrangian $\mathcal{L}:Z\times X^{\prime}\to \mathbb{R}$ defined by
\[
\begin{aligned}
\mathcal{L}(y, \lambda):=& \frac{1}{2}a(y,y) + b(y,\lambda)- l(y).
\end{aligned}
\]
\item We have the estimate 
\begin{equation}
\Vert y\Vert_Z= \Vert y\Vert_{L^2(q_T)}\leq \Vert y_{obs}\Vert_{L^2(q_T)} ,  \quad \Vert \lambda\Vert_{X^{\prime}} \leq 2\sqrt{C_{\Omega,T}+\eta}  \Vert y_{obs}\Vert_{L^2(q_T)}. \label{estimate_lambda} 
\end{equation}
\end{enumerate}
\end{theorem}

\textsc{Proof-} We use classical results for saddle point problems  (see \cite{brezzi_new}, chapter 4). 

We easily obtain the continuity of the bilinear form $a$ over $Z\times Z$, the continuity of bilinear $b$ over $Z\times X^{\prime}$ and the continuity of the linear form $l$ over $Z$. In particular, we get 
\begin{equation}
\Vert l\Vert_{Z^{\prime}}= \Vert y_{obs}\Vert_{L^2(q_T)}, \qquad \Vert a\Vert_{(Z\times Z)^{\prime}}=1, \quad \Vert b\Vert_{(Z\times X^{\prime})^{\prime}} = \eta^{-1/2}. \label{simple_ineq}
\end{equation}
Moreover, the kernel $\mathcal{N}(b)=\{y\in Z;\  b(y,\lambda)=0 \quad \forall \lambda\in X^{\prime}\}$ coincides with $W$: we easily get 
\begin{equation}
a(y,y)=\Vert y \Vert^2_Z, \quad \forall y\in \mathcal{N}(b)=W.  \nonumber
\end{equation}
Therefore, in view of \cite[Theorem 4.2.2]{brezzi_new}, it remains to check the inf-sup constant property : $\exists \delta>0$ such that
\begin{equation}
\inf_{\lambda\in X^{\prime}} \sup_{y \in Z}  \frac{b(y,\lambda)}{\Vert y \Vert_Z \Vert \lambda\Vert_{X^{\prime}}} \geq \delta.
\end{equation}
We proceed as follows. For any fixed $\lambda\in X^{\prime}$, we define $y$ as the unique solution of 
\begin{equation}
Ly=-\Delta\lambda \,\,\,\textrm{in} \,\,\, Q_T,  \quad (y(\cdot,0),y_t(\cdot,0))=(0,0)\,\,\, \textrm{on}\,\,\, \Omega ,\quad y=0\,\,\, \textrm{on}\,\,\, \Sigma_T.   \label{infsup_y0}
\end{equation}  
We get $b(y,\lambda)=\Vert \lambda\Vert^2_{X^{\prime}}$
and
\begin{equation}
\Vert y\Vert^2_Z= \Vert y\Vert^2_{L^2(q_T)} + \eta \Vert \lambda \Vert_{X^{\prime}}^2.   \nonumber
\end{equation}
Using (\ref{continuity_constant}), the estimate $\Vert y\Vert_{L^2(q_T)} \leq \sqrt{C_{\Omega,T}} \Vert \lambda \Vert_{X^{\prime}}$  
 implies that $y\in Z$ and that 
\[
\sup_{y\in Z}  \frac{b(y,\lambda)}{\Vert y \Vert_Z \Vert \lambda \Vert_{X^{\prime}}} \geq \frac{1}{\sqrt{C_{\Omega,T}+\eta}}>0
\]
leading to the result with $\delta=(C_{\Omega,T}+\eta)^{-1/2}$. 

The third point is the consequence of classical estimates  (see \cite{brezzi_new}, Theorem 4.2.3.) : 
\begin{equation}
\Vert y\Vert_{Z}\leq \frac{1}{\alpha_0} \Vert l\Vert_{Z^{\prime}}, \quad \Vert \lambda \Vert_{X^{\prime}} \leq \frac{1}{\delta}\biggl(1+\frac{\Vert a \Vert}{\alpha_0} \biggr) \Vert l\Vert_{Z^{\prime}} \nonumber
\end{equation}
where 
\begin{equation}
\alpha_0:=\inf_{y\in \mathcal{N}(b)}  \frac{a(y,y)}{\Vert y\Vert^2_Z}.  \label{eq:alpha0}
\end{equation}
Estimates (\ref{simple_ineq}) and the equality $\alpha_0=1$ lead to the results.  Eventually, from (\ref{simple_ineq}), we obtain that 
\begin{equation}
\Vert \lambda\Vert_{X^{\prime}} \leq \frac{2}{\delta} \Vert y_{obs}\Vert_{L^2(q_T)}   \nonumber
\end{equation}
and that $\delta \geq (C_{\Omega,T}+\eta)^{-1/2}$ to get (\ref{estimate_lambda}). \Fin

In practice, it is very convenient to "augment" the Lagrangian (see \cite{fortinglowinski}) and consider instead the Lagrangian $\mathcal{L}_r$ defined for any $r>0$ by 
\begin{equation}
\begin{aligned}
& \mathcal{L}_r(y,\lambda):=\frac{1}{2}a_r(y,y)+b(y,\lambda)-l(y), \\
& a_r(y,y):=a(y,y)+r\Vert Ly\Vert^2_{X}. 
\end{aligned}
 \nonumber
\end{equation}
Since $a_r(y,y)=a(y,y)$ on $W$, the Lagrangian $\mathcal{L}$ and $\mathcal{L}_r$ share the same saddle-point. The positive number $r$ is an augmentation parameter. 

\begin{remark}\label{rk_lambda_sys}
Assuming additional hypotheses on the regularity of the solution $\lambda$, precisely $L\lambda\in L^2(Q_T)$ and $(\lambda,\lambda_t)_{\vert t=0,T} \in H^1_0(\Omega)\times L^2(\Omega)$,
we easily prove, writing the optimality condition for $\mathcal{L}$, that the multiplier $\lambda$ satisfies the following relations :
%
%
\begin{equation}
\label{system_lambda}
\left\{
\begin{aligned}
& L\lambda=-(y-y_{obs})\, 1_\omega\quad \textrm{in}\quad Q_T , \quad  \lambda=0\quad \textrm{in}\quad \Sigma_T, \\
& \lambda=\lambda_{t}=0\quad\textrm{on}\ \Omega \times \{0,T\}.
\end{aligned}
\right.
\end{equation}
Therefore, $\lambda$ (defined in the weak sense) is an exact controlled solution of the wave equation through the control $-(y-y_{obs})\, 1_{\omega}\in L^2(q_T)$. 
\begin{itemize}
\item  If $y_{obs}$ is the restriction to $q_T$ of a solution of (\ref{eq:wave}), then the unique multiplier $\lambda$ must vanish almost everywhere. In that case, 
we have  $\sup_{\lambda\in\Lambda}\inf_{y\in Y} \mathcal{L}_r(y,\lambda) = \inf_{y\in Y} \mathcal{L}_r(y,0)=\inf_{y\in Y} J_r(y)$ with 
\begin{equation}
\label{def_Jyr}
J_r(y):=\frac{1}{2}\Vert y-y_{obs}\Vert^2_{L^2(Q_T)} + \frac{r}{2}\Vert Ly\Vert^2_X.
\end{equation}
The corresponding variational formulation is then : find $y\in Z$ such that 
\begin{equation}
a_r(y,\overline{y})=\jjntqT y\,\overline{y} \, dxdt + r\int_0^T  \langle Ly, \, L\overline{y} \rangle_{H^{-1}(\Omega)} \,dt = l(\overline{y}), \quad \forall \overline{y}\in Z.  \nonumber
\end{equation}
\item In the general case, the mixed formulation can be rewritten as follows: find $(z,\lambda)\in Z\times X^{\prime}$ solution of 
\begin{equation}
\left\{
\begin{aligned}
\langle P_r y, P_r \overline{y}\rangle_{X\times L^2(q_T)} + \langle L\overline{y},\lambda \rangle_{X,X^{\prime}}   & = \langle (0,y_{obs}),P_r \overline{y} \rangle_{X\times L^2(q_T)}, \quad \forall \overline{y}\in Z, \\
\langle L\overline{y},\lambda \rangle_{X,X^{\prime}}   & = 0, \quad \forall \lambda\in X^{\prime}
\end{aligned}
\right.
\end{equation}
with $P_r y:=(\sqrt{r} L\,y,y_{\vert q_T})$. This approach may be seen as generalization of the \eqref{FV} problem (see (\ref{FV})), where the variable $\lambda$ is adjusted automatically (while the choice of the parameter $\eps$ in (\ref{FV}) is in general a delicate issue).
\end{itemize}

\end{remark}

\

System (\ref{system_lambda}) can be used to define a equivalent saddle-point formulation, very suitable at the numerical level.  Precisely, we introduce - in view of (\ref{system_lambda}) - the space $\Lambda$ by 
$$
\nonumber
\begin{aligned}
\Lambda:=\{\lambda: \lambda\in C([0,T]; & H_0^1(\Omega))\cap C^1([0,T]; L^2(\Omega)), \\
&L\lambda\in L^2(Q_T), \lambda(\cdot,0)=\lambda_t(\cdot,0)=0\}.
\end{aligned}
$$
Endowed with the scalar product  $\langle \lambda,\overline{\lambda}\rangle_{\Lambda}:= \jjntQT (\lambda\, \overline{\lambda}+L\lambda L\overline{\lambda})  \, dxdt$, we check that $\Lambda$ is a Hilbert space.  Then, for any parameter $\alpha\in (0,1)$, we consider the following mixed formulation : find $(y,\lambda)\in Z\times \Lambda$ such that 

\begin{equation} \label{eq:mfalpha}
\left\{
\begin{array}{rcll}
\noalign{\smallskip} a_{r,\alpha}(y, \overline{y}) + b_{\alpha}(\overline{y}, \lambda) & = & l_{1,\alpha}(\overline{y}), & \qquad \forall \overline{y} \in Z \\
\noalign{\smallskip} b_{\alpha}(y, \overline{\lambda}) - c_{\alpha}(\lambda,\overline{\lambda})& = & l_{2,\alpha}(\overline{\lambda}), & \qquad \forall \overline{\lambda} \in \Lambda,
\end{array}
\right.
\end{equation}
where
\begin{align}
& a_{r,\alpha} : Z \times Z \to \mathbb{R},  \quad a_{r,\alpha}(y,\overline{y}) := (1-\alpha)\jjntqT y\overline{y}\,dxdt + r \int_0^T  (Ly, L\overline{y})_{H^{-1}(\Omega)}dt, \nonumber\\
& b_{\alpha}: Z\times \Lambda   \to \mathbb{R},  \quad b_{\alpha}(y,\lambda) := \int_0^T  \langle \lambda,Ly \rangle_{H^1_0(\Omega),H^{-1}(\Omega)} dt-\alpha \jjntqT y\, L\lambda\,dxdt, \nonumber\\  
& c_{\alpha}: \Lambda \times \Lambda  \to \mathbb{R},  \quad c_{\alpha}(\lambda, \overline{\lambda}) := \alpha\jjntQT L\lambda\,L\overline{\lambda}, \, dxdt \nonumber\\
& l_{1,\alpha}: Z\to \mathbb{R},  \quad l_{1,\alpha}(y) := (1-\alpha)\jjntqT y_{obs} \, y\, dxdt, \nonumber\\
& l_{2,\alpha}: \Lambda\to \mathbb{R},  \quad l_{2,\alpha}(\lambda) := -\alpha \jjntqT  y_{obs}\, L \lambda\,dxdt.  \nonumber
\end{align}

From the symmetry of $a_{r,\alpha}$ and $c_{\alpha}$, we easily check that this formulation corresponds to the saddle point problem : 
\begin{equation}
\nonumber
\left\{
\begin{aligned}
& \sup_{\lambda\in \Lambda} \inf_{y\in Z} \mathcal{L}_{r,\alpha}(y,\lambda), \\
& \mathcal{L}_{r,\alpha}(y,\lambda):=\mathcal{L}_r(y,\lambda) - \frac{\alpha}{2} \Vert L\lambda + (y-y_{obs})1_{\om}\Vert_{L^2(Q_T)}^2. 
\end{aligned}
\right.
\end{equation}
\begin{proposition}
Under the hypothesis \eqref{iobs}, for any $\alpha\in (0,1)$, the formulation (\ref{eq:mfalpha}) is well-posed. Moreover, the unique pair $(y,\lambda)$ in $Z\times \Lambda$ satisfies 
\begin{equation}
\theta_1 \Vert y\Vert_Z^2 + \theta_2 \Vert \lambda\Vert_{\Lambda}^2 \leq  \biggl(\frac{(1-\alpha)^2}{\theta_1} + \frac{\alpha^2}{\theta_2}\biggr)\Vert y_{obs}\Vert^2_{L^2(q_T)}.  \label{estimate_solalpha}
\end{equation}
with 
\begin{equation}
\theta_1:=\min\biggl(1-\alpha, \frac{r}{\eta}\biggr) , \quad \theta_2:=\frac{1}{2}\min\biggl(\alpha, \frac{1}{C_{\Omega,T}}\biggr).   \nonumber
\end{equation}
 \end{proposition}
 \textsc{Proof-} We easily get the continuity of the bilinear forms $a_{r,\alpha}$, $b_{\alpha}$ and $c_{\alpha}$: 
 $$
 \begin{aligned}
& \vert a_{r,\alpha}(y,\overline{y})\vert \leq \max(1-\alpha,\frac{r}{\eta}) \Vert y\Vert_Z \Vert\overline{y}\Vert_Z, \quad \forall y,\overline{y}\in Z, \\
& \vert b_{\alpha}(y,\lambda)\vert  \leq \max(\alpha,\frac{1}{\sqrt{\eta}}) \Vert y\Vert_{Z} \Vert \lambda\Vert_{\Lambda}, \quad \forall y\in Z, \forall \lambda\in \Lambda  ,\\
& \vert c_{\alpha}(\lambda,\overline{\lambda}) \leq \alpha \Vert \lambda\Vert_{\Lambda} \Vert \overline{\lambda}\Vert_{\Lambda}, \quad \forall \lambda, \overline{\lambda}\in \Lambda 
 \end{aligned}
 $$
 and of the linear form $l_1$ and $l_2$ : $\Vert l_1\Vert_{Z^{\prime}}=(1-\alpha)\Vert y_{obs}\Vert_{L^2(q_T)}$ and 
$\Vert l_2\Vert_{\Lambda^{\prime}}=\alpha\Vert y_{obs}\Vert_{L^2(q_T)}$.   
 
 Moreover, since $\alpha\in (0,1)$, we also obtain the coercivity of $a_{r,\alpha}$ and of $c_{\alpha}$: precisely, 
$$
\begin{aligned}
& a_{r,\alpha}(y,y)\geq \min\biggl(1-\alpha, \frac{r}{\eta}\biggr) \Vert y\Vert^2_Z, \quad\forall y\in Z, \\
& c_{\alpha}(\lambda,\lambda)\geq \min\biggl(\alpha m, \frac{1-m}{C_{\Omega,T}}\biggr) \Vert \lambda\Vert^2_{\Lambda} \quad\forall \lambda\in \Lambda, \quad \forall m\in (0,1).
\end{aligned}
$$
The result \cite[Prop 4.3.1]{brezzi_new}  implies the well-posedness and the estimate (\ref{estimate_solalpha}) taking $m=1/2$. \Fin

The $\alpha$-term in $\mathcal{L}_{r,\alpha}$ is a stabilization term: it ensures a coercivity property of $\mathcal{L}_{r,\alpha}$ with respect to the variable $\lambda$ and automatically the well-posedness. In particular, there is no need to prove any inf-sup property for the application $b_{\alpha}$.

\begin{proposition}
If the solution $(y,\lambda)\in Z\times X^\prime$ of (\ref{eq:mf}) enjoys the property $\lambda\in \Lambda$, then the solutions of (\ref{eq:mf}) and (\ref{eq:mfalpha}) coincide. 
\end{proposition}
\textsc{Proof-} The hypothesis of regularity and the relation (\ref{system_lambda}) imply that the solution $(y,\lambda)\in Z\times X^\prime$ of (\ref{eq:mf}) is also a solution of (\ref{eq:mfalpha}). The result then follows from the uniqueness of the two formulations. \Fin

\subsection{Dual formulation of the extremal problem (\ref{eq:mf})} \label{sec2_dual}

As discussed at length in \cite{NC-AM-mixedwave}, we may also associate to the extremal problem \eqref{P} an equivalent problem involving only the variable $\lambda$. Again, this is particularly interesting at the numerical level. This requires a strictly positive augmentation parameter $r$.

For any $r>0$, let us define the linear operator $\mathcal{P}_r$ from $X^{\prime}$ into $X^{\prime}$ by 
\[
\mathcal{P}_r\lambda:= -\Delta^{-1} (Ly), \quad \forall \lambda\in X^{\prime}
\]
where $y \in Z$ is the unique solution to
\begin{equation}\label{eq:imageA}
a_r(y, \overline y) = b(\overline y, \lambda), \quad \forall \overline y \in Z.   
\end{equation}
The assumption $r>0$ is necessary here in order to guarantee the well-posedness of (\ref{eq:imageA}). Precisely, for any $r>0$, the form $a_r$ defines a norm equivalent to the norm on $Z$.

The following important lemma holds:
\begin{lemma}\label{propA}
For any $r>0$, the operator $\mathcal{P}_r$ is a strongly elliptic, symmetric isomorphism from $X^{\prime}$ into $X^{\prime}$.
\end{lemma}
\textsc{Proof-} From the definition of $a_r$, we easily get that $\Vert \mathcal{P}_r\lambda\Vert_{X^{\prime}}\leq r^{-1} \Vert \lambda\Vert_{X^\prime}$ and the continuity of $\mathcal{P}_r$. Next, consider any $\lambda^{\prime}\in X^{\prime}$ and denote by $y^{\prime}$ the corresponding unique solution of (\ref{eq:imageA}) so that $\mathcal{P}_r\lambda^{\prime}:=-\Delta^{-1}(Ly^{\prime})$. Relation (\ref{eq:imageA}) with $\overline{y}=y^{\prime}$ then implies that 
\begin{equation}
\int_0^T \langle\mathcal{P}_r\lambda^{\prime},\lambda\rangle_{H_0^1(\Omega)} \,dt = a_r(y,y^{\prime})    \label{arAlambda}
\end{equation}
and therefore the symmetry and positivity of $\mathcal{P}_r$. The last relation with $\lambda^{\prime}=\lambda$ and the observability estimate (\ref{iobs}) imply that $\mathcal{P}_r$ is also positive definite.

Finally, let us check the strong ellipticity of $\mathcal{P}_r$, equivalently that the bilinear functional $(\lambda,\lambda^{\prime})\to \int_0^T \langle\mathcal{P}_r\lambda,\lambda^{\prime}\rangle_{H^1_0(\Omega),H^1_0(\Omega)}\,dt$ is $X^{\prime}$-elliptic. Thus we want to show that 
\begin{equation}\label{ellipticity_A}
\int_0^T \langle\mathcal{P}_r\lambda,\lambda\rangle_{H_0^1(\Omega)}\,dt \geq C \Vert \lambda\Vert^2_{X^{\prime}}, \quad\forall \lambda\in X^{\prime}
\end{equation}
for some positive constant $C$. Suppose that (\ref{ellipticity_A}) does not hold; there exists then a sequence $\{\lambda_n\}_{n\geq 0}$ of $X^{\prime}$ such that 

\[
\Vert \lambda_n\Vert_{X^{\prime}}=1, \quad\forall n\geq 0, \qquad \lim_{n\to\infty} \int_0^T \langle\mathcal{P}_r\lambda_n,\lambda_n\rangle_{H_0^1(\Omega)}\,dt=0.
\]
Let us denote by $y_n$ the solution of (\ref{eq:imageA}) corresponding to $\lambda_n$. From (\ref{arAlambda}), we then obtain that 
\begin{equation}
\lim_{n\to\infty} r \Vert Ly_n\Vert^2_X+ \Vert y_n\Vert^2_{L^2(q_T)}=0 .\label{limit}
\end{equation}
From (\ref{eq:imageA}) with $y=y_n$ and $\lambda=\lambda_n$, we have
\begin{equation} \label{phinlambdan}
\int_0^T  \left\langle r (-\Delta^{-1}) L y_n-\lambda_n , (-\Delta^{-1})L\overline{y} \right\rangle_{H^1_0(\Omega)}\,dt + \jjntqT y_n \overline{y} dx\,dt =0, \quad \forall \overline{y} \in Z.
\end{equation}
We define the sequence $\{\overline{y}_n\}_{n\geq 0}$ as follows : 
\begin{equation}\nonumber
\left\{
\begin{aligned}
& L\overline{y}_n = r\,Ly_n+\Delta^{-1}\lambda_n, & &\textrm{in}\quad Q_T, \\
& \overline{y}_n=0, & &\textrm{in}\quad \Sigma_T,\\
& \overline{y}_n(\cdot,0)=\overline{y}_{n,t}(\cdot,0)=0, & &\textrm{in}\quad \Omega,
\end{aligned}
\right.
\end{equation}
so that, for all $n$, $\overline{y}_n$ is the solution of the wave equation with zero initial data and source term $r Ly_n +\Delta\lambda_n$ in $X$. Using again (\ref{continuity_constant}), we get $\Vert \overline{y}_n\Vert_{L^2(q_T)} \leq \sqrt{C_{\Omega,T}} \Vert r Ly_n+\Delta\lambda_n\Vert_X$, so that $\overline{y}_n\in Z$. Then, using (\ref{phinlambdan}) with $\overline{y}=\overline{y}_n$ we get 
\[
\Vert r(-\Delta^{-1})Ly_n-\lambda_n\Vert_{X^{\prime}} \leq \sqrt{C_{\Omega,T}} \Vert y_n\Vert_{L^2(q_T)}.
\]
Then, from (\ref{limit}), we conclude that $\lim_{n\to +\infty} \Vert \lambda_n\Vert_{X^{\prime}}=0$ leading to a contradiction and to the strong ellipticity of the operator $\mathcal{P}_r$. \Fin

The introduction of the operator $\mathcal{P}_r$ is motivated by the following proposition~: 

\begin{proposition}\label{prop_equiv_dual}
For any $r>0$, let $y_0\in Z$ be the unique solution of 
\[
a_r(y_0,\overline{y})= l(\overline{y}), \quad \forall \overline{y}\in Z
\]
and let $J_r^{\star\star}:X^{\prime}\to X^{\prime}$ be the functional defined by 
\[
J_r^{\star\star}(\lambda) = \frac{1}{2} \int_0^T \langle\mathcal{P}_r \lambda, \lambda\rangle_{H_0^1(\Omega)}dt - b(y_0, \lambda).
\]
The following equality holds : 
\begin{equation}\nonumber
\sup_{\lambda\in X^{\prime}}\inf_{y\in Z} \mathcal{L}_r(y,\lambda) = - \inf_{\lambda\in X^{\prime}} J_r^{\star\star}(\lambda)\quad + \mathcal{L}_r(y_0,0).
\end{equation}
\end{proposition}
The proof is classical and we refer for instance to \cite{NC-AM-mixedwave} in a similar context.
This proposition reduces the search of $y$, solution of problem \eqref{P}, to the minimization of $J_r^{\star\star}$. The well-posedness is a consequence of the ellipticity of the operator $\mathcal{P}_r$. 


\begin{remark}
The results of this section apply if the distributed observation on $q_T$ is replaced by a Neumann boundary observation on a sufficiently large subset $\Sigma_T$ of 
$\partial\Omega\times (0,T)$ (i.e.  assuming $\frac{\partial y}{\partial \nu}=y_{obs}\in L^2(\Sigma_T)$ is known on $\Sigma_T$). This is due to the following generalized observability inequality: there exists a positive constant $C_{obs}=C(\omega,T,\Vert c\Vert_{C^1(\overline{\Omega})},\Vert d\Vert_{L^{\infty}(\Omega)})$ such that the following estimate holds : 
\begin{equation}
\Vert y(\cdot,0),y_t(\cdot,0)\Vert^2_{H_0^1(\Omega)\times L^2(\Omega)} \leq C_{obs} \biggl( \biggl\Vert \frac{\partial y}{\partial \nu}\biggr\Vert^2_{L^2(\Sigma_T)}+ \Vert Ly\Vert^2_{L^2(Q_T)}  \biggr), \quad \forall y\in Z  
\end{equation}
which holds if the triplet $(Q_T,\Sigma_T,T)$ satisfies the geometric condition as before (we refer to \cite{NC-EFC-AM} and the references therein). Actually, it suffices to re-define the form 
$a$ in \eqref{eq:a} by $a(y,\overline{y}):=\int\!\!\!\int_{\Sigma_T} \frac{\partial y}{\partial\nu}  \frac{\partial \overline{y}}{\partial\nu} \, d\sigma dx$ and the form $l$ by $l(y):=\int\!\!\!\int_{\Sigma_T} \frac{\partial y}{\partial\nu} y_{obs} \, d\sigma dx$ for all $y,\overline{y}\in Z$.
\end{remark}

\begin{remark}
We emphasize that the mixed formulation (\ref{eq:mf}) has a structure very closed to the one we get when we address - using the same approach - the null controllability of (\ref{eq:wave}): more precisely, the control of minimal $L^2(q_T)$-norm which drives to rest the initial data $(y_0,y_1)\in H_0^1(\Omega)\times L^2(\Omega)$ is given by $v=\ph\, 1_{q_T}$ where $(\ph,\lambda)\in \Phi\times L^2(0,T; H_0^1(\Omega))$ solves the mixed formulation
\begin{equation} \label{eq:mfi}
\left\{
\begin{array}{rcll}
\noalign{\smallskip} a(\ph, \overline{\ph}) + b(\overline{\ph}, \lambda) & = & l(\overline{\ph}), & \qquad \forall \overline{\ph} \in \Phi \\
\noalign{\smallskip} b(\ph, \overline{\lambda}) & = & 0, & \qquad \forall \overline{\lambda} \in L^2(0,T; H^1_0(\Omega)),
\end{array}
\right.
\end{equation}
where
\begin{align}
& a : \Phi \times \Phi \to \mathbb{R},  \quad a(\ph, \overline{\ph}) = \jjntqT \ph(x, t) \overline{\ph}(x, t) \, dx \, dt \nonumber\\
& b: \Phi \times L^{2}(0,T; H^1_0(\Omega))  \to \mathbb{R},  \quad b(\ph, \lambda) = \int_0^T \langle L \ph, \lambda \rangle_{H^{-1},H^1_0} dt \nonumber\\
& l: \Phi \to \mathbb{R},  \quad l(\ph) = -\langle\ph_t(\cdot, 0),y_0\rangle_{H^{-1}(\Omega),H_0^1(\Omega)} + \int_0^1 \ph(\cdot, 0)\, y_1 \,dx. \nonumber
\end{align}
with $\Phi = \left\{\ph \in L^2(Q_T), \, \ph = 0 \text{ on } \Sigma_T \text{ such that } L\ph \in L^2(0,T; H^{-1}(\Omega))\right\}$. We refer to \cite{NC-AM-mixedwave}.
\end{remark}

\begin{remark}
Reversing the order of priority between the constraint  $y-y_{obs}=0$ in $L^2(q_T)$ and $Ly-f=0$ in $X$, a possibility could be to minimize the functional 
$y\to \Vert Ly-f\Vert_{X}$ over $y\in Z$ subject to the constraint $y-y_{obs}=0$ in $L^2(q_T)$ via the introduction of a Lagrange multiplier in $L^2(q_T)$. 
The proof of the following inf-sup property : there exists $\delta>0$ such that 
$$
\inf_{\lambda\in L^2(q_T)} \sup_{y\in Z} \frac{\int\!\!\!\int_{q_T} \lambda y \,dxdt }{\Vert \lambda\Vert_{L^2(q_T)} \Vert y\Vert_Y}   \geq \delta
$$
associated to the corresponding mixed-formulation is however unclear. 
If a $\eps$-term is added as in (\ref{FV}), this difficulty disappears (we refer again to the book \cite{Klibanov-book}).
\end{remark}


\section[Recovering the source and the solution from a partial observation]{Recovering the source and the solution from a partial observation: a mixed re-formulation of the problem}\label{recovering_y_f}

Given a partial observation $y_{obs}$ of the solution on the subset $q_T\subset Q_T$, we now consider the reconstruction of the full solution as well as the source term $f$ assumed in $X$. We assume that the initial data $(y_0,y_1)\in \boldsymbol{H}$ are unknown. 

The situation is different with respect to the previous section, since without additional assumption on $f$, the couple $(y,f)$ is not unique. Consider the case of a source $f$ supported in a set which is near $\partial\Omega\times (0,T)$  and disjoint from $q_T$: from the finite propagation of the solution, the source $f$ will not affect the solution $y$ in $q_T$. 
On the other hand, the determination of a couple $(y,f)$ which solves (\ref{eq:wave}) such that $y$ coincides with $y_{obs}$ is straightforward : it suffices to "extend" $y$ on $Q_T\setminus q_T$ appropriately to preserve the boundary conditions, then compute $Ly$ and recover a source term. However, we emphasize that, on a practical viewpoint, the extension of $y_{obs}$ out of $q_T$ is not obvious. Moreover, this strategy does not offer any control on the object $f$. 

We briefly show that we can apply the method developed in Section \ref{recovering_y} which allows a robust reconstruction and then consider the case of uniqueness via additional condition on $f$.

We assume again that \eqref{iobs} holds. We note $Y:=Z\times X$ and define on $Y$ the bilinear form, for any $\eps,\eta>0$
\begin{equation}
\label{eq:pseta2}
\begin{aligned}
\langle (y,f),(\overline{y},\overline{f}) \rangle_Y:= & \jjntqT y\,\overline{y}\,dxdt + \eta\int_0^T \langle Ly-f,L\overline{y}-f\rangle_{H^{-1}(\Omega)}dt \\
&+ \eps \int_0^T \langle f, \overline{f} \rangle_{H^{-1}(\Omega)} dt, \qquad \forall (y,f), (\overline{y},\overline{f})\in Y.
\end{aligned}
\end{equation}
In view of (\ref{iobs}), this bilinear form defines a scalar product over $Y$. Moreover, endowed to this scalar product, we easily obtain that $Y$ is a Hilbert space (we refer to \cite{CC-NC-AM}). We note the corresponding norm by $\Vert (y,f)\Vert_{Y}:=\sqrt{((y,f),(y,f))_Y}$.

Then, for any $\eps>0$, we consider the following extremal problem : 
\begin{equation}
\nonumber
(\mathcal{P}_{\eps})\quad 
\left\{
\begin{aligned}
& \inf  J_{\eps}(y,f):=  \frac{1}{2}\Vert y-y_{obs}\Vert^2_{L^2(q_T)} + \frac{\eps}{2}\Vert f\Vert^2_X, \\
& \textrm{subject to}\quad (y,f)\in W
\end{aligned}
\right.
\end{equation}
where $W$ is the closed subspace of $Y$ defined by 
$
W:=\{(y,f)\in Y; \,  Ly-f=0 \,\, \textrm{in}\,\, X\}
$
and endowed with the norm of $Y$: precisely, it follows that 
$$\Vert (y,f)\Vert_{W}:=\sqrt{\Vert y\Vert_{L^2(q_T)}^2 +\eps \Vert f\Vert^2_X}, \qquad \forall (y,f)\in W.$$

The extremal problem $(\mathcal{P}_{\eps})$ is well posed : the functional $J_{\eps}$ is continuous over $W$, is strictly convex and is such that $J_{\eps}(y,f)\to +\infty$ as $\Vert (y,f)\Vert_W\to \infty$. Note also that the solution of $(\mathcal{P}_{\eps})$ in $W$, depends on $\eps$ but not on $\eta$. 

Remark also that if $\eps=0$, then $J_{\eps}$ is \textit{a priori} only convex leading possibly to distinct minima. This justifies the introduction of the $\eps$-term in the functional $J_{\eps}$. We emphasize however that the $\eps$-term is not a regularization term as it does not improve the regularity of the state $y$. 

Eventually, from (\ref{iobs}), the solution $(y_\eps,f_\eps)$ in $W$ of $(\mathcal{P}_{\eps})$ satisfies $(y_\eps(\cdot,0),y_{\eps,t}(\cdot,0))\in \boldsymbol{H}$, so that problem ($\mathcal{P}_{\eps}$) is again equivalent to the minimisation of $J_{\eps}$ with respect to $(y_0,y_1,f)\in \boldsymbol{H}\times X$. 


Proceeding as in Section \ref{recovering_y}, we introduce a Lagrangian multiplier $\lambda_{\eps}\in X^{\prime}$ and the following mixed formulation: find $((y_{\eps},f_{\eps}), \lambda_{\eps})\in Y\times X^{\prime}$ solution of 
\begin{equation} \label{eq:mfeps}
\left\{
\begin{array}{rcll}
\noalign{\smallskip} a_{\eps}((y_{\eps},f_{\eps}), (\overline{y},\overline{f})) + b((\overline{y},\overline{f}), \lambda_{\eps}) & = & l(\overline{y},\overline{f}), & \qquad \forall (\overline{y},\overline{f}) \in Y \\
\noalign{\smallskip} b((y_{\eps},f_{\eps}), \overline{\lambda}) & = & 0, & \qquad \forall \overline{\lambda} \in X^{\prime},
\end{array}
\right.
\end{equation}
where
\begin{align}
\label{eq:aeps} & a_{\eps} : Y \times Y \to \mathbb{R},  \quad a_{\eps}((y,f), (\overline{y},\overline{f})) := \jjntqT y\overline{y}\,dxdt + \eps (f,\overline{f})_X,
\\
\label{eq:bf} & b: Y \times X^{\prime}  \to \mathbb{R},  \quad b((y,f),\lambda) := \int_0^T \langle \lambda, Ly-f\rangle_{H_0^1(\Omega),H^{-1}(\Omega)} \, dt,\\
\label{eq:lf} & l: Y \to \mathbb{R},  \quad l(y,f) := \jjntqT y_{obs} \,y\, dxdt. 
\end{align}
\begin{theorem}\label{th:mfeps} Under the hypothesis \eqref{iobs}, the following hold :
\begin{enumerate}
\item The mixed formulation (\ref{eq:mfeps}) is well-posed.
\item The unique solution $((y_{\eps},f_{\eps}), \lambda_{\eps}) \in Y\times X^{\prime}$ is the saddle-point of the Lagrangian $\mathcal{L}_{\eps}:Y\times X^\prime\to \mathbb{R}$ defined by
\[
\mathcal{L}_{\eps}((y,f), \lambda):= \frac{1}{2}a_{\eps}((y,f),(y,f)) + b((y,f),\lambda)- l(y,f).
\]
Moreover, the pair $(y_{\eps},f_{\eps})$ solves the extremal problem $(\mathcal{P}_{\eps})$.

\item The following estimates hold :  
\begin{equation}
\Vert (y_{\eps},f_{\eps}) \Vert_Y= \left(\Vert y_{\eps}\Vert^2_{L^2(q_T)}+\eps  \Vert f_{\eps}\Vert^2_{} \right)^{1/2}\leq \Vert y_{obs}\Vert_{L^2(q_T)}  \label{ineqy_th1}
\end{equation}
and 
\begin{equation}
\Vert \lambda_{\eps}\Vert_{L^2(Q_T)} \leq 2\sqrt{C_{\Omega,T}+\eta}  \Vert y_{obs}\Vert_{L^2(q_T)}  \label{estimate_lambdaeps}
\end{equation}
for some constant $C_{\Omega,T}>0$.
\end{enumerate}
\end{theorem}

The proof is very closed to the proof of Theorem \ref{th:mf}.  In particular, the obtention of the inf-sup property is obtained by taking, for any $\lambda\in X^{\prime}$, $f=0$ and $y$ as in (\ref{infsup_y0}) so that the inf-sup constant
\begin{equation}
\delta_{\eps}:=\inf_{\lambda\in X^{\prime}} \sup_{(y,f)\in Y}  \frac{b((y,f),\lambda)}{\Vert (y,f) \Vert_Y \Vert \lambda\Vert_{X^{\prime}}}
\end{equation}
is bounded by above by $(C_{\Omega,T}+\eta)^{-1/2}$ uniformly with respect to $\eps$.

Remark in particular that the inequality (\ref{ineqy_th1}) implies that, at the optimality, since $\eps>0$, the equality $\Vert y-y_{obs}\Vert_{L^2(q_T)}=0$ can not hold if $f_{\eps}\neq 0$.


\begin{remark}\label{inf-sup_proof2}
We may also prove the inf-sup property using the variable $f$: for any $\lambda\in X^{\prime}$, we set $y=0$ and  $f=\Delta\lambda\in X$. We get 
\begin{equation}
\nonumber
\sup_{(y,f)\in Y} \frac{b((y,f),\lambda)}{\Vert (y,f) \Vert_Y \Vert \lambda \Vert_{X^{\prime}} } \geq   \frac{b((0,\Delta \lambda),\lambda)}{\Vert (0,\Delta \lambda) \Vert_{Y} \Vert \lambda \Vert_{X^{\prime}} }= \frac{1}{\sqrt{\eps+\eta}}
\end{equation}
so that $\delta_{\eps}\geq (\eps+\eta)^{-1/2}$. Therefore, the estimate 
\begin{equation}
\Vert \lambda_{\eps}\Vert_{X^{\prime}} \leq \frac{2}{\delta_{\eps}} \Vert y_{obs}\Vert_{L^2(q_T)}   \nonumber
\end{equation}
implies that 
\begin{equation}
\Vert \lambda_{\eps}\Vert_{X^{\prime}} \leq 2\sqrt{\eps+\eta}\Vert y_{obs}\Vert_{L^2(q_T)}.\label{estimate_lambdaepsnew}
\end{equation}
This argument is valid if and only if $f$ is distributed everywhere in $Q_T$. \Fin
\end{remark}

\begin{remark}\label{FVeps}
The estimate (\ref{estimate_lambdaepsnew}) implies that the multiplier $\lambda_{\eps}$ vanishes in $X^{\prime}$ as $\eps+\eta\to 0^+$ (recall that $\eps$ and $\eta$ can be chosen arbitrarily small in (\ref{eq:pseta})). 
\end{remark}

\begin{remark}\label{rq_stab_eps}
\begin{itemize}
\item[]
\item[(a)] Assuming enough regularity on the solution $\lambda_{\eps}$, precisely that $L\lambda_\eps\in L^2(Q_T)$ and $(\lambda,\lambda_t)_{t=0,T}\in H_0^1(\Omega)\times L^2(\Omega)$, we easily check that the multiplier $\lambda_{\eps}$ satisfies the following relations :
\begin{equation}
\nonumber
\left\{
\begin{aligned}
& L\lambda_{\eps}=-(y_{\eps}-y_{obs})_{1_\omega}, \quad Ly_{\eps}-f_{\eps}=0, \quad  \eps f_{\eps}+\Delta \lambda_{\eps}=0\quad\textrm{in}\quad Q_T,  \\
& \lambda_{\eps}=0\quad \textrm{in}\quad \Sigma_T, \\
& \lambda_{\eps}=\lambda_{\eps,t}=0\quad\textrm{on}\ \Omega \times \{0,T\}.
\end{aligned}
\right.
\end{equation}
Therefore, $\lambda_{\eps}$ is an exact controlled solution of the wave equation through the control $-(y_{\eps}-y_{obs})\, 1_{\omega}$ and from (\ref{estimate_lambdaepsnew})
implies that
\begin{equation}
\Vert y_{\eps}-y_{obs}\Vert_{L^2(q_T)} \to 0 \quad \textrm{as}\quad  \eps\to 0^+.  \label{lim_eps}
\end{equation}
Remark however that $f_{\eps}$ may not be bounded in $X^{\prime}$ uniformly w.r.t. $\eps$ (contrarily to the sequence $(\sqrt{\eps} f_{\eps})_{\eps>0}$).

\item[(b)] The equality $Ly_{\eps}=f_{\eps}$ becomes $\eps Ly_{\eps}=-\Delta\lambda_{\eps}$ and leads to $L(\eps \Delta^{-1}Ly_{\eps})=-L\lambda_{\eps}=(y_{\eps}-y_{obs})_{1_\omega}$. Finally, 
$y_{\eps}$ solves, at least in $\mathcal{D}^{\prime}$, the boundary value problem
\begin{equation}
\nonumber
\left\{
\begin{aligned}
& L(\eps (-\Delta ^{-1})L y_{\eps} ) + y_{\eps}\,1_{\omega} = y_{obs}\, 1_{\omega}, \quad \textrm{in}\quad Q_T, \\
& (\eps Ly_{\eps})=(\eps Ly_{\eps})_t=0,  \quad \textrm{in}\quad \Omega\times \{0,T\} \\
& y_{\eps} =0, \quad \textrm{on}\quad \Sigma_T 
\end{aligned}
\right.
\end{equation}
or equivalently the variational formulation: find $y_{\eps}\in Z$ (see \eqref{refZ}) solution of 
\begin{equation}
\eps \int_0^T \langle Ly_{\eps},L\overline{y} \rangle_{H^{-1}(\Omega)} dt + \jjntqT y_{\eps} \, \overline{y} \,dxdt = \jjntqT y_{obs} \overline{y}\, dxdt, \quad\forall \overline{y}\in Z  \label{FV_eps}
\end{equation}
which actually can be obtained directly from the cost $J_{\eps}$, replacing from the beginning $f$ by the term $Ly$. From the Lax-Milgram lemma, (\ref{FV_eps}) is well-posed and the following estimates hold :  
$$
\Vert y_{\eps}\Vert_{L^2(q_T)}\leq \Vert y_{obs}\Vert_{L^2(q_T)}, \quad \sqrt{\eps}\Vert Ly_{\eps}\Vert_{X}\leq \Vert y_{obs}\Vert_{L^2(q_T)}.
$$ 
This kind of variational formulation involving the fourth order term $Ly_{\eps}Ly$ has been derived and used in \cite{NC-EFC-AM} in a controllability context.
\end{itemize}
\end{remark}

For any $\eps>0$ and any $y_{obs}\in L^2(q_T)$, the method allows to recover a couple $(y_{\eps},f_{\eps})$ such that $L y_{\eps}=f_{\eps}$ in $Q_T$ and $y_{\eps}$ is closed to $y_{obs}$ (see (\ref{lim_eps})). In view of the loss of uniqueness, we have no information on the limit of the sequence as $\eps\to 0$: the sequence may be unbounded at the limit in $L^2(Q_T)\times L^2(Q_T)$ even if $y_{obs}$ is the restriction to $q_T$ of a solution of (\ref{eq:wave}).

\begin{remark} Contrarily to the inf-sup property, the coercivity of $a_{\eps}$ over $\mathcal{N}(b)$ does not hold uniformly with respect to $\eps$. Recall that the $\eps$-term has been introduced to get a norm for $Y$. This enforces us to add this term in the mixed formulation.
\end{remark}

\begin{remark}
\textit{A fortiori}, if the initial condition $(y_0,y_1)\in \boldsymbol{H}$ is known, one may recover the pair $(y,f)\in Y$ from $y_{obs}$ and $(y_0,y_1)$. The procedure is similar; it suffices to define two additional Lagrange multipliers $(\lambda_1,\lambda_2)\in L^2(\Omega)\times H^1_0(\Omega)$ to deal with the constraint $y(\cdot,0)=y_0$ and $y_t(\cdot,0)=y_1$ respectively. The extremal problem is now :
\begin{equation}
\inf_{(y,f)\in W}   J_{\eps}(y,f):=  \frac{1}{2}\Vert y-y_{obs}\Vert^2_{L^2(q_T)} + \frac{\eps}{2}\Vert f\Vert^2_{X^{\prime}}   \nonumber
\end{equation}
where $W$ is the closed subspace of $Y$ defined by 
$$
W:=\{(y,f)\in Y; \,  Ly-f=0 \,\, \textrm{in}\,\, X^{\prime}, (y(\cdot,0),y_t(\cdot,0))=(y_0,y_1)\,\, \textrm{in}\,\, \boldsymbol{H}\}. 
$$
The corresponding mixed formulation is :
find $((y_{\eps},f_{\eps}),  (\lambda_{\eps},\lambda_{\eps,1},\lambda_{\eps,2})) \in Y\times \Lambda$ solution of 
\begin{equation} \label{eq:mfeps1}
\left\{
\begin{array}{rcll}
\noalign{\smallskip} a_{\eps}((y_\eps,f_\eps), (\overline{y},\overline{f})) + b((\overline{y},\overline{f}), (\lambda_\eps,\lambda_{\eps,1},\lambda_{\eps,2})) & = & l_1(\overline{y},\overline{f}), & \qquad \forall (\overline{y},\overline{f}) \in Y \\
\noalign{\smallskip} b((y_\eps,f_\eps), (\overline{\lambda},\overline{\lambda_1},\overline{\lambda_2})) & = & l_2(\overline{\lambda},\overline{\lambda_1},\overline{\lambda_2}), & \qquad \forall (\overline{\lambda},\overline{\lambda_1},\overline{\lambda_2}) \in \Lambda,
\end{array}
\right.
\end{equation}
where $a_\varepsilon$ is given by \eqref{eq:aeps} and
\begin{align}
 & b: Y\times \Lambda  \to \mathbb{R},  \quad b((y,f),(\lambda,\lambda_1,\lambda_2)) := \jjntQT \lambda (Ly-f) \, dx dt  \nonumber\\
& \hspace{5cm} + \langle y(\cdot,0),\lambda_1\rangle_{L^2(\Omega)} + \langle y_t(\cdot,0),\lambda_2 \rangle_{H^{-1}(\Omega),H_0^1(\Omega)} \nonumber\\
& l_1: Y \to \mathbb{R},  \quad l_1(y,f) := \jjntqT y_{obs} \overline{y}\, dxdt \nonumber\\
& l_2: \Lambda \to \mathbb{R},  \quad l_2(\lambda,\lambda_1,\lambda_2) := \langle y_0,\lambda_1\rangle_{L^2(\Omega)} + \langle y_1,\lambda_2\rangle_{H^{-1}(\Omega),H_0^1(\Omega)} \nonumber
\end{align}
with $\Lambda:=X^{\prime}\times L^2(\Omega)\times H^1_0(\Omega)$.
Using the estimate (\ref{iobs}), we easily show that this formulation is well-posed.  \Fin
\end{remark}

In view of Remark \ref{rq_stab_eps} (a), we may also associate to the mixed formulation \eqref{eq:mfeps} a stabilized version, similarly to (\ref{eq:mfalpha}). 

Again, it is very convenient to "augment" the Lagrangian (see \cite{fortinglowinski}) and consider instead the Lagrangian $\mathcal{L}_{\eps,r}$ defined for any $r>0$ by 
\begin{equation}
\begin{aligned}
& \mathcal{L}_{\eps,r}((y,f),\lambda):=\frac{1}{2}a_{\eps,r}((y,f),(y,f))+b(y,\lambda)-l(y,f), \\
& a_{\eps,r}((y,f),(y,f)):=a_{\eps}((y,f),(y,f))+r\Vert Ly-f\Vert^2_X. 
\end{aligned}
 \nonumber
\end{equation}
Since $a_{\eps}(y,y)=a_{\eps,r}(y,y)$ on $W$, the Lagrangian $\mathcal{L}_{\eps}$ and $\mathcal{L}_{\eps,r}$ share the same saddle-point. The positive number $r$ is an augmentation parameter. Similarly, proceeding as in Section \ref{sec2_dual}, we may also associate to the saddle-point problem $\sup_{\lambda\in X^{\prime}}\inf_{(y,f)\in Y} \mathcal{L}_{r,\eps}((y,f),\lambda)$ a dual problem, which again reduces the search of the couple $(y_{\eps},f_{\eps})$, solution of problem $(\mathcal{P}_{\eps})$, to the minimization of a elliptic functional in $\lambda_{\eps}$. 

\begin{proposition}\label{prop_equiv_dual_eps}
For any $r>0$, let $(y_0,f_0)\in Y$ be the unique solution of 
\[
a_{\eps,r}((y_0,f_0),(\overline{y},\overline{f}))= l(\overline{y},\overline{f}), \quad \forall (\overline{y},\overline{f})\in Y
\]
and let $\mathcal{P}_{\eps,r}$ be the strongly elliptic and symmetric operator from $X^{\prime}$ into $X^{\prime}$ defined by 
$\mathcal{P}_{\eps,r}\lambda:= -\Delta^{-1}(Ly-f)$ where $(y,f)\in Y$ is the unique solution to 
\begin{equation}\label{eq:imageAeps}
a_{\eps,r}((y,f), (\overline y,\overline f)) = b((\overline y,\overline f), \lambda), \quad \forall (\overline y,\overline f) \in Y.   
\end{equation}
Then, the following equality holds
\begin{equation}\nonumber
\sup_{\lambda\in X^{\prime}}\inf_{(y,f)\in Y} \mathcal{L}_{\eps,r}((y,f),\lambda) = - \inf_{\lambda\in X^{\prime}} J_{\eps,r}^{\star\star}(\lambda)\quad + \mathcal{L}_{\eps,r}((y_0,f_0),0).
\end{equation}
where $J_{\eps,r}^{\star\star}:X^{\prime}\to X^{\prime}$ is the functional defined by 
\[
J_{\eps,r}^{\star\star}(\lambda) = \frac{1}{2} \int_0^T (\mathcal{P}_{\eps,r} \lambda, \lambda)_{H_0^1(\Omega)} \,dt - b((y_0,f_0), \lambda).
\]
\end{proposition}

Compared to the previous section, the additional unknown $f_{\epsilon}$ on the problem guarantees that the term $\Vert y_{\eps}-y_{obs}\Vert_{L^2(q_T)}$
vanishes at the limit in $\eps$, for any $y_{obs}\in L^2(q_T)$, be a restriction of a solution of $(\ref{eq:wave})$ or not. The situation is different if additional assumption on $f$ enforces the uniqueness of the pair $(y,f)$ (we refer to \cite{yamamoto} and the references therein).


\section{Numerical Analysis of the mixed formulations} \label{sec_numer}

\subsection{Numerical approximation of the mixed formulation (\ref{eq:mf})}

We consider the numerical analysis of the mixed formulation (\ref{eq:mf}), assuming $r>0$. We follow \cite{NC-AM-mixedwave}, to which we refer for the details. 

Let $Z_h$ and $\Lambda_h$ be two finite dimensional spaces parametrized by the variable $h$ such that $Z_h\subset Z, \Lambda_h\subset X^{\prime}$ for every $h>0$.
Then, we can introduce the following approximated problems: find the $(y_h,\lambda_h)\in Z_h\times \Lambda_h$ solution of 
\begin{equation} \label{eq:mfh}
\left\{
\begin{array}{rcll}
\noalign{\smallskip} a_r(y_h,\overline{y}_h) + b(\overline{y}_h, \lambda_h) & = & l(\overline{y}_h), & \qquad \forall \overline{y}_h\in Z_h \\
\noalign{\smallskip} b(y_h, \overline{\lambda}_h) & = & 0, & \qquad \forall \overline{\lambda}_h \in \Lambda_h.
\end{array}
\right.
\end{equation}
The well-posedness of this mixed formulation is again a consequence of two properties: the coercivity of the bilinear form $a_r$ on the subset 
\[\mathcal{N}_h(b)=\{y_h\in Z_h; b(y_h,\lambda_h)=0\,\qquad \forall \lambda_h\in \Lambda_h\}. 
\]
 Actually, from the relation $a_r(y,y)\geq (r/\eta)\Vert y\Vert_Z^2$ for all $y\in Z$, the form $a_r$ is coercive on the full space $Z$, and so a fortiori on $\mathcal{N}_h(b)\subset Z_h\subset Z$. The second property is a discrete inf-sup condition. We note $\delta_h>0$ by
\begin{equation}
\delta_h:=\inf_{\lambda_h\in \Lambda_h}\sup_{y_h\in Z_h}  \frac{b(y_h,\lambda_h)}{\Vert \lambda_h\Vert_{X^{\prime}}\Vert y_h\Vert_Z}. \label{infsupdiscret}
\end{equation}
For any fixed $h$, the spaces $Z_h$ and $\Lambda_h$ are of finite dimension so that the infimum and supremum in (\ref{infsupdiscret}) are reached: moreover, from the property of the bilinear form $a_r$, it is standard to check that $\delta_h$ is strictly positive. Consequently, for any fixed $h>0$, there exists a unique couple $(y_h,\lambda_h)$ solution of (\ref{eq:mfh}). 

We then have the following estimate.

\begin{proposition} \label{prop_estimateh_eps0}
Let $h>0$. Let $(y,\lambda)$ and $(y_h,\lambda_h)$ be the solution of (\ref{eq:mf}) and of (\ref{eq:mfh}) respectively. Let $\delta_h$ the discrete inf-sup constant defined by (\ref{infsupdiscret}). Then,
\begin{align}
\label{eq:est_e01}
& \Vert y-y_h \Vert_{Z} \leq  2\biggl(1+\frac{1}{\sqrt{\eta} \delta_h}\biggr)d(y,Z_h)+\frac{1}{\sqrt{\eta}} d(\lambda,\Lambda_h), \\
\label{eq:est_e02}
& \Vert \lambda-\lambda_h \Vert_{X^{\prime}}  \leq  \biggl(2+\frac{1}{\sqrt{\eta}\delta_h}\biggr)\frac{1}{\delta_h}d(y,Z_h)+
\frac{3}{\sqrt{\eta} \delta_h} d(\lambda,\Lambda_h)
\end{align}
where  
$ d(\lambda,\Lambda_h):=\inf_{\lambda_h\in \Lambda_h} \Vert \lambda-\lambda_h\Vert_{X^{\prime}}$
and
$$
\begin{aligned}
d(y,Z_h):=& \inf_{y_h\in Z_h}\Vert y-y_h\Vert_Z\\
=& \inf_{y_h \in Z_h}   \biggl( \Vert y-y_h\Vert^2_{L^2(q_T)}+  \eta \Vert L(y-y_h) \Vert^2_X \biggr)^{1/2}.
\end{aligned}
$$
\end{proposition}
\textsc{Proof-} From the classical theory of approximation of saddle point problems (see \cite[Theorem 5.2.2]{brezzi_new}) we have that
\begin{align}
\nonumber
\Vert y-y_h \Vert_{Z} \leq & \left( \frac{2 \| a_r \|_{(Z\times Z)^{\prime}}}{\alpha_0} + \frac{2 \| a_r \|_{(Z\times Z)^{\prime}}^\frac12 \| b \|_{(Z\times X^\prime)^{\prime}}}{\alpha_0^\frac12 \delta_h} \right) d(y,Z_h) \\
& +\frac{\Vert b\Vert_{(Z\times X^\prime)^{\prime}}}{\alpha_0} d(\lambda,\Lambda_h)
\end{align}
and 
\begin{align}\nonumber
\Vert \lambda-\lambda_h \Vert_{X^{\prime}} \leq  &
\left( \frac{2\| a_r \|_{(Z\times Z)^{\prime}}^\frac32}{\alpha_0^\frac12 \delta_h} + \frac{\| a_r \|_{(Z\times Z)^{\prime}} \| b \|_{(Z\times X^\prime)^{\prime}}}{\delta_h^2} \right) d(y,Z_h) \\
& + \frac{3 \| a_r \|^\frac12 \| b \|_{(Z\times X^\prime)^{\prime}}}{\alpha_0^\frac12 \delta_h} d(\lambda,\Lambda_h).
\end{align}
Since, $\Vert a_r\Vert_{(Z\times Z)^{\prime}}=\alpha_0=1$; $\Vert b\Vert_{(Z\times \Lambda)^{\prime}}= \frac{1}{\sqrt{\eta}}$, the result follows.  \Fin 
\begin{remark}
For $r=0$, the discrete mixed formulation (\ref{eq:mfh}) is not well-posed over $Z_h\times \Lambda_h$ because the form $a_{r=0}$ is not coercive over the discrete kernel of $b$: the equality $b(y_h,\lambda_h)=0$ for all $\lambda_h\in \Lambda_h$ does not imply in general that $L y_h$ vanishes. Therefore, the term $r \Vert Ly_h\Vert^2_{X}$, which appears in the Lagrangian $\mathcal{L}_r$, may be understood as a stabilization term: for any $h>0$, it ensures the uniform coercivity of the form $a_r$ and vanishes at the limit in $h$. We also emphasize that this term is not a regularization term as it does not add any regularity on the solution $y_h$. 
\end{remark}


Let $n_h=\dim Z_h, m_h=\dim \Lambda_h$ and let the real matrices $A_{r,h}\in \mathbb{R}^{n_h,n_h}$, $B_{h}\in \mathbb{R}^{m_h,n_h}$, $J_h\in \mathbb{R}^{m_h,m_h}$ and $L_h\in \mathbb{R}^{n_h}$ be defined by   
\begin{equation}
\label{def_matrix}
\left\{
\begin{aligned}
& a_r( y_h, \overline{y_h})=   \langle A_{r,h} \{y_h\}, \{\overline{y_h}\} \rangle_{\mathbb{R}^{n_h},\mathbb{R}^{n_h}} & \forall y_h,\overline{y_h}\in Z_h,\\
& b(y_h,\lambda_h)=  \langle B_{h} \{y_h\}, \{\lambda_h\} \rangle_{\mathbb{R}^{m_h},\mathbb{R}^{m_h}} & \forall y_h\in Z_h, \lambda_h\in \Lambda_h,\\
& \jjntQT \lambda_h\overline{\lambda_h}\,dx\,dt=  \langle J_h \{\lambda_h\}, \{\overline{\lambda_h}\} \rangle_{\mathbb{R}^{m_h},\mathbb{R}^{m_h}} & \forall \lambda_h,\overline{\lambda_h}\in \Lambda_h, \\
& l(y_h)= \langle L_h,\{y_h\} \rangle_{\mathbb{R}^{n_h}} & \forall y_h\in Z_h,
\end{aligned}
\right.
\end{equation}
where $\{y_h\}\in \mathbb{R}^{n_h}$ denotes the vector associated to $y_h$ and $\langle \cdot,\cdot \rangle_{\mathbb{R}^{n_h},\mathbb{R}^{n_h}}$ the usual scalar product over $\mathbb{R}^{n_h}$. With these notations, the problem (\ref{eq:mfh}) reads as follows: find $\{y_h\}\in \mathbb{R}^{n_h}$ and $\{\lambda_h\}\in \mathbb{R}^{m_h}$ such that 
\begin{equation} \label{matrixmfh}
\left(
\begin{array}{cc}
A_{r,h} &  B_h^T \\
B_h & 0   
\end{array}
\right)_{\mathbb{R}^{n_h+m_h,n_h+m_h}}
\left(
\begin{array}{c}
\{y_h\}  \\
\{\lambda_h\}    
\end{array}
\right)_{\mathbb{R}^{n_h+m_h}}  =
\left(
\begin{array}{c}
L_h \\
0   
\end{array}
\right)_{\mathbb{R}^{n_h+m_h}}.
\end{equation}
The matrix $A_{r,h}$ as well as the mass matrix $J_h$ are symmetric and positive definite for any $h>0$ and any $r>0$. On the other hand, the matrix of order $m_h+n_h$ in (\ref{matrixmfh}) is symmetric but not positive definite. We use exact integration methods developed in \cite{dunavant} for the evaluation of the coefficients of the matrices. The system (\ref{matrixmfh}) is solved using the direct LU decomposition method. 

\subsubsection{$C^1$-finite elements and order of convergence for $N=1$}

The finite dimensional and conformal space $Z_h$ must be chosen such that $Ly_h$ belongs to $X=L^2(0,T; H^{-1}(\Omega))$ for any $y_h\in Z_h$. This is guaranteed, for instance, as soon as $\ph_h$ possesses second-order derivatives in $L^2_{loc}(Q_T)$. As in \cite{NC-AM-mixedwave}, we consider a conformal approximation based on functions continuously differentiable with respect to both variables $x$ and $t$. 

We introduce a triangulation $\mathcal{T}_h$ such that $\overline{Q_T}=\cup_{K\in \mathcal{T}_h} K$ and we assume that $\{\mathcal{T}_h\}_{h>0}$ is a regular family. We note 
$
h:=\max\{\textrm{diam}(K), K\in \mathcal{T}_h\}, 
$
where $\textrm{diam}(K)$ denotes the diameter of $K$. Then, we introduce the space $Z_h$ as follows : 
\begin{equation}
Z_h=\{y_h\in Z \subset C^1(\overline{Q_T}):  z_h\vert_K\in \mathbb{P}(K) \quad \forall K\in \mathcal{T}_h, \,\, z_h=0 \,\,\textrm{on}\,\,\Sigma_T\},   \label{defZh}
\end{equation}
where $\mathbb{P}(K)$ denotes an appropriate space of functions in $x$ and $t$. In this work, we consider two choices, in the one-dimensional setting (for which $\Omega\subset \mathbb{R}$, $Q_T\subset \mathbb{R}^2$): 
\begin{enumerate}
\item The \textit{Bogner-Fox-Schmit} (BFS for short) $C^1$-element defined for rectangles. It involves $16$ degrees of freedom, namely the values of $y_h, y_{h,x}, y_{h,t}, y_{h,xt}$ on the four vertices of each rectangle $K$. 
Therefore $\mathbb{P}(K)= \mathbb{P}_{3,x}\otimes \mathbb{P}_{3,t}$ where $\mathbb{P}_{r,\xi}$ is by definition the space of polynomial functions of order $r$ in the variable $\xi$. We refer to \cite[ch. II, sec. 9, p. 94]{ciarletfem}.
\item The reduced \textit{Hsieh-Clough-Tocher} (HCT for short) $C^1$-element defined for triangles. This is a so-called composite finite element and involves $9$ degrees of freedom, namely, the values of $y_h,y_{h,x}, y_{h,t}$ on the three vertices of each triangle $K$. We refer to \cite[ch. VII, sec. 46, p. 285]{ciarletfem} and to \cite{bernadouHCT,meyer} where the implementation is discussed.
\end{enumerate}
We also define the finite dimensional space
\[
\Lambda_h=\{\lambda_h\in C^0(\overline{Q_T}),  \lambda_h\vert_K\in \mathbb{Q}(K) \quad \forall K\in \mathcal{T}_h \}.
\]
where $\mathbb{Q}(K)$ denotes the space of affine functions both in $x$ and $t$ on the element $K$.

For any $h>0$, we have $Z_h\subset Z$ and $\Lambda_h\subset X^{\prime}$.

We then have the following result: 
\begin{proposition}[BFS element for $N=1$ - Rate of convergence for the norm $Z$]   \label{propeps0}
Let $h>0$, let $k\in \{1,2\}$ be a positive integer. Let $(y,\lambda)$ and $(y_h,\lambda_h)$ be the solution of (\ref{eq:mf}) and (\ref{eq:mfh}) respectively. If the solution $(y,\lambda)$ belongs to $H^{k+2}(Q_T)\times H^k(Q_T)$, then there exists two positives constants 
\[ K_i=K_i(\Vert y\Vert_{H^{k+2}(Q_T)}, \Vert c \Vert_{C^1(\overline{Q_T})},\Vert d \Vert_{L^{\infty}(Q_T)}), \qquad i \in \{ 1,2\},\]
independent of $h$, such that 
\begin{align}
\Vert y-y_h \Vert_Z & \leq  K_1 \frac{h^{k-1}}{\sqrt{\eta}}\biggl(  (\sqrt{\eta}+\frac{1}{\delta_h})  (h^3+\sqrt{\eta} h) + 1 \biggr), \label{eq:yyh0}\\
 \Vert \lambda-\lambda_h \Vert_{X^{\prime}} &  \leq  K_2 \frac{h^{k-1}}{\sqrt{\eta}\delta_h}\biggl(  (\sqrt{\eta}+\frac{1}{\delta_h})  (h^3+\sqrt{n}h) + 1  \biggr)\label{eq:llh0}.
\end{align}
\end{proposition}
\par\noindent
\textsc{Proof -}
 From \cite[ch. III, sec. 17]{ciarletfem}, for any $\lambda\in H^k(Q_T)$, $k\leq 2$, there exists $C_1=C_1(\Vert \lambda\Vert_{H^k(Q_T)})$
such that 
\begin{equation} \label{eq:est0h}
\Vert \lambda-\Pi_{\Lambda_h,\mathcal{T}_h}(\lambda)\Vert_{X^{\prime}} \leq  C_1 h^{k-1}, \quad \forall h>0
\end{equation}
where $\Pi_{\Lambda_h,\mathcal{T}_h}$ designates the interpolant operator from $X^{\prime}$ to $\Lambda_h$ associated to the regular mesh $\mathcal{T}_h$. 
Similarly, for any $y\in H^{k+2}(Q_T)$, there exist $C_2=C_2(\Vert y\Vert_{H^{k+2}(Q_T)})$ such that for every $h > 0$ we have
\begin{equation} \label{eq:est1h}
\Vert y-\Pi_{Z_h,\mathcal{T}_h}(y)\Vert_{L^2(Q_T)} \leq  C_2  h^{k+2}, \quad  \Vert y-\Pi_{Z_h,\mathcal{T}_h}(y)\Vert_{H^2(Q_T)} \leq  C_2  h^k. 
\end{equation}
Then, observing that 
\begin{equation} \label{eq:est2h}
\Vert Ly-Ly_h\Vert_X \leq K(\Vert c\Vert_{C^1(\overline{Q_T})}, \Vert d\Vert_{L^{\infty}(Q_T)})\Vert y-y_h\Vert_{H^2(Q_T)},
\end{equation}
for some positive constant $K$, we get that 
\begin{equation}
\begin{aligned}
d(y,Z_h)&=\inf_{y_h \in Z_h} \left( \|y-y_h\|^2_{L^2(q_T)} + \eta \|Ly - Ly_h\|^2_{X}\right)^2 \\
&\leq  C_2 \biggl((h^{k+2}) ^2 +  \eta K^2 (h^k)^2\biggr)^{1/2} \\
& \leq C_2  (h^{k+2}+\sqrt{n}K\,h^k)
\end{aligned}
\label{eq:ch}
\end{equation}
and then from Proposition \ref{prop_estimateh_eps0}, we get that 
\begin{equation}
 \Vert y-y_h \Vert_{Z} \leq  2\biggl(1+\frac{1}{\sqrt{\eta} \delta_h}\biggr)C_2(h^{k+2}+\sqrt{n}K\,h^k)+\frac{1}{\sqrt{\eta}} C_1h^{k-1}.
\end{equation}
Similarly, 
 \begin{equation*}
\Vert \lambda-\lambda_h \Vert_{X^{\prime}}\leq  \biggl(2+\frac{1}{\sqrt{\eta} \delta_h}\biggr)\frac{1}{\delta_h} C_2 (h^{k+2}+\sqrt{n}K\,h^k)+\frac{3}{\sqrt{\eta}\delta_h} C_1h^{k-1}.
\end{equation*}
From the last two estimates, we obtain the conclusion of the proposition. 
\Fin

It remains now to deduce the convergence of the approximated solution $y_h$ for the $L^2(Q_T)$ norm: this is done using the observability estimate (\ref{iobs}). Precisely, we write that  $(y-y_h)$ solves 
\begin{equation}
\nonumber
\left\{
\begin{aligned}
& L(y-y_h)=-L y_h \quad \text{in } Q_T \\
& ((y-y_h), (y-y_h)_t)(0) \in \boldsymbol{H} \\
& y-y_h=0  \quad \text{on } \Sigma_T.
\end{aligned}
\right.
\end{equation}
Therefore using (\ref{estimate_qTQT}), there exists a constant $C(C_{\Omega,T},C_{obs})$ such that 
\begin{equation}
\Vert y-y_{h} \Vert^2_{L^2(Q_T)} \leq C(C_{\Omega,T},C_{obs}) ( \Vert y-y_h \Vert^2_{L^2(q_T)} +  \Vert L y_h\Vert^2_{X})  \nonumber
\end{equation}
from which we deduce, in view of the definition of the norm $Y$, that 
\begin{equation}
\Vert y-y_h \Vert_{L^2(Q_T)} \leq C(C_{\Omega,T},C_{obs})  \max(1,\frac{2}{\sqrt{\eta}}) \Vert y-y_h\Vert_Z.   \label{eq:estyL2}
\end{equation}
Eventually, by coupling (\ref{eq:estyL2}) and Proposition \ref{propeps0}, we obtain the following result :

\begin{theorem}[BFS element for $N=1$ - Rate of convergence for the norm $L^2(Q_T)$]   \label{theps0}
Assume that the hypothesis \eqref{iobs} holds. Let $h>0$, let $k\in \{1,2\}$ be a positive integer and let $\eta < 1$. Let $(y,\lambda)$ and $(y_h,\lambda_h)$ be the solution of (\ref{eq:mf}) and (\ref{eq:mfh}) respectively. If the solution $(y,\lambda)$ belongs to $H^{k+2}(Q_T)\times H^k(Q_T)$, then there exists two positives constant $K=K(\Vert y\Vert_{H^{k+2}(Q_T)}, \Vert c \Vert_{C^1(\overline{Q_T})},$ $\Vert d \Vert_{L^{\infty}(Q_T)},C_{\Omega,T},C_{obs})$, independent of $h$, such that 
\begin{equation}
\Vert y-y_h \Vert_{L^2(Q_T)}  \leq  K \max(1,\frac{2}{\sqrt{\eta}}) \frac{h^{k-1}}{\sqrt{\eta}}\biggl(  (\sqrt{\eta}+\frac{1}{\delta_h})  (h^3+\sqrt{\eta} h) + 1 \biggr).   \label{estimate_L2_eps0}
\end{equation}
\end{theorem}

\begin{remark}\label{better_estimate}
Estimate (\ref{estimate_L2_eps0}) is not fully satisfactory as it depends on the constant $\delta_h$. 
In view of the complexity of both the constraint $Ly=0$ and of the structure of the space $Z_h$, the theoretical estimation of the constant $\delta_h$ with respect to $h$ is a difficult problem. However, as discussed at length in  \cite[Section 2.1]{NC-AM-mixedwave}, $\delta_h$ can be evaluated numerically for any $h$,  as the solution of the following generalized eigenvalue problem (taking $\eta=r$, so that $a_r(y,y)$ is exactly $\Vert y \Vert^2_{Z}$):   
\begin{equation}\label{eigenvalue}
\delta_h = \inf\biggl\{\sqrt{\delta}:   B_h A_{r,h}^{-1} B_h^T  \{\lambda_h\} = \delta \,J_h \{\lambda_h\}, \quad \forall\, \{\lambda_h\}\in \mathbb{R}^{m_h}\setminus\{0\}\biggr\}
\end{equation}
where the matrix $A_{r,h}$, $B_h$ and $J_h$ are defined in (\ref{def_matrix}).

Table \ref{tab:infsup} reports the values of $\delta_h$ for $r=1$ and $r=h^{-2}$ for several values of $h$, $T=2$, $\omega=(0.1,0.3)$ and the BFS element. 
As in \cite{NC-AM-mixedwave} where the boundary situation is considered with more details, these values suggests that, asymptotically in $h$, the constant $\delta_{r,h}$ behaves like :
\begin{equation}
\delta_{r,h}\approx C_r \frac{1}{\sqrt{r}} \quad \textrm{as} \quad h \to 0^+   \label{behavior_deltah_eps0}
\end{equation}
with $C_r>0$, a uniformly bounded constant w.r.t. $h$. 
\begin{table}[http]
\centering
\begin{tabular}{|c|cccc|}
\hline
$h$  &  $7.01\times 10^{-2}$  & $3.53\times 10^{-2}$ & $1.76\times 10^{-2}$ & $8.83\times 10^{-3}$  \tabularnewline

\hline 

$r=1$ & $3.58$ & $3.48$ & $3.42$ & $3.40$ \tabularnewline

$r=h^{-2}$ & $2.53\times 10^{-1}$ & $1.23\times 10^{-1}$ & $6.05\times 10^{-2}$ & $3.01\times 10^{-2}$ \tabularnewline

\hline
\end{tabular}
\caption{$\eps=0$:  $T=2$ - $\delta_{r,h}$ for $r=1$ and $r=h^{-2}$ with respect to $h$. }
\label{tab:infsup}
\end{table}

Consequently, in view of \ref{behavior_deltah_eps0}, the right hand side of the estimate (\ref{estimate_L2_eps0}) of $\Vert y-y_{h}\Vert_{L^2(Q_T)}$ behaves, taking $\eta=r$ and $r>1$
so that $\max(1,\frac{1}{\sqrt{r}})=1$, like 
$$
\Vert y-y_{h}\Vert_{L^2(Q_T)} \leq K h^{k-1} \biggl(\sqrt{r}h+\frac{1}{\sqrt{r}}\biggr)
$$
and reaches its minimum for $r=1/h$, leading to $\Vert y-y_{h}\Vert_{L^2(Q_T)} \leq K h^{k-1/2}$.
\end{remark}

Eventually, when the space $Z_h$ is based on the HCT element, Theorem \ref{theps0} and Remark \ref{better_estimate} still hold for $k=1$. From \cite[ch. VII, sec. 48, p. 295]{ciarletfem}, we use that, for $k\in \{0,1\}$, there exists a constant $C_2>0$ such
\begin{equation}
\Vert y-\Pi_{Z_h,\mathcal{T}_h}(y)\Vert_{L^2(Q_T)} \leq  C_2  h^{k+2}, \quad  \Vert y-\Pi_{Z_h,\mathcal{T}_h}(y)\Vert_{H^2(Q_T)} \leq  C_2  h^k.
\end{equation}
Then, we use that the error $\Vert y-y_h\Vert_{L^2(Q_T)}$ is again controlled by the error on the Lagrange multiplier $\lambda$ through the term $d(\lambda,\Lambda_h)$ in (\ref{eq:est_e01}) to conclude.

\subsection{Numerical approximation of the mixed formulation (\ref{eq:mfalpha})}

We address the numerical approximation of the \textit{stabilized} mixed formulation (\ref{eq:mfalpha}) with $\alpha \in (0,1)$ and $r>0$. Let $h$ be a real parameter. 
Let $Z_h$ and $\widetilde\Lambda_h$ be two finite dimensional spaces such that 
$$
Z_h\subset Z, \quad \widetilde\Lambda_h\subset \Lambda, \qquad \forall h>0.
$$
The problem (\ref{eq:mfalpha}) becomes : find $(y_{h}, \lambda_{h}) \in Z_h\times \widetilde\Lambda_h$ solution of 
\begin{equation} \label{eq:mfalphah}
\left\{
\begin{array}{rcll}
\noalign{\smallskip} a_{r,\alpha}(y_h, \overline{y}_h) + b_{\alpha}( \lambda_h, \overline{y}_h) & = & l_{1,\alpha}(\overline{y}_h), & \qquad \forall \overline{y}_h \in Z_h \\
\noalign{\smallskip} b_{\alpha}( \overline{\lambda}_h, y_h) - c_{\alpha}(\lambda_h,\overline{\lambda}_h)& = & l_{2,\alpha}(\overline{\lambda}_h), & \qquad \forall \overline{\lambda}_h \in \widetilde\Lambda_h,
\end{array}
\right.
\end{equation}

Proceeding as in the proof of \cite[Theorem 5.5.2]{brezzi_new}, we first easily show that the following estimate holds .
\begin{lemma} \label{lemmastab} Let $(y,\lambda)\in Y\times \Lambda$ be the solution of (\ref{eq:mfalpha}) and $(y_h,\lambda_h)\in Z_h\times \widetilde{\Lambda_h}$ be the solution of (\ref{eq:mfalphah}). Then we have, 
\begin{align}
\frac{1}{4}\theta_1 \Vert y-y_h\Vert_{Z}^2 +   \frac{1}{4}\theta_2 \Vert \lambda-\lambda_h\Vert_{\widetilde\Lambda}^2  \leq  & \biggl(\frac{\Vert a_{r,\alpha}\Vert^2}{\alpha_a}+\frac{\Vert b_\alpha\Vert^2}{\alpha_c}+\frac{\theta_1}{2}\biggr) \inf_{\overline{y}_h\in Z_h} \Vert\overline{y}_h-y\Vert_Z^2 \nonumber \\
& + \biggl(\frac{\Vert b_{\alpha}\Vert^2}{\theta_1}+\frac{\alpha^2}{\theta_2}+\frac{\theta_2}{2}\biggr) \inf_{\overline{\lambda}_h\in \widetilde{\Lambda_h}}\Vert \overline{\lambda}_h-\lambda\Vert_{\Lambda}^2
\end{align}
with  $\Vert a_{r,\alpha}\Vert\leq \max(1-\alpha, \eta^{-1}r)$,  $\Vert b_{\alpha}\Vert\leq \max(\eta^{-1/2},\alpha)$. Parameters $\theta_1$ and $\theta_2$ are defined in (\ref{estimate_solalpha}).
\end{lemma}
Concerning the space $\widetilde\Lambda_h$, since $L\lambda_h$ should belong to $L^2(Q_T)$, a natural choice is 
\begin{equation}
\widetilde\Lambda_h=\{\lambda\in Z_h; \lambda(\cdot,0)=\lambda_t(\cdot,0)=0\}.  \label{choice_lambdah_alpha}
\end{equation}
where $Z_h\subset Z$ is defined by (\ref{defZh}).
Then, using Lemma \ref{lemmastab} and the estimate (\ref{eq:ch}), we obtain the following result. 
\begin{proposition}[BFS element for $N=1$ - Rate of convergence - Stabilized formulation]
Let $h>0$, let $k \leq 2$ be a positive integer and let $\alpha\in (0,1)$. Let $(y,\lambda)$ and $(y_h,\lambda_h)$ be the solution of (\ref{eq:mfalpha}) and (\ref{eq:mfalphah}) respectively. If $(y,\lambda)$ belongs to $H^{k+2}(Q_T)\times H^{k+2}(Q_T)$, then there exists a positive constant $K=K(\Vert y\Vert_{H^{k+2}(Q_T)}, \Vert c \Vert_{C^1(\overline{Q_T})},\Vert d \Vert_{L^{\infty}(Q_T)},\alpha,r,\eta)$  independent of $h$, such that 
\begin{equation}
\begin{aligned}
\Vert y-y_h \Vert_Z + \Vert \lambda-\lambda_h \Vert_{\Lambda} \leq K h^k.
\end{aligned}
\end{equation}
\end{proposition}

In particular, arguing as in the previous section, we get 

\begin{theorem}[Rate of convergence for the norm $L^2(Q_T)$ Stabilized formulation] 
Assume that the hypothesis \eqref{iobs} holds. Let $h>0$, let an integer $k\leq 2$. Let $(y,\lambda)$ and $(y_{h},\lambda_{h})$ be the solution of (\ref{eq:mfalpha}) and (\ref{eq:mfalphah}) respectively. If the solution $(y, \lambda)$ belongs to $H^{k+2}(Q_T)\times H^{k+2}(Q_T)$, then there exist a positive constant $K=K(\Vert y\Vert_{H^{k+2}(Q_T)}, \Vert \lambda\Vert_{H^{k+2}(Q_T)}, \Vert c \Vert_{C^1(\overline{Q_T})},\Vert d \Vert_{L^{\infty}(Q_T)},\alpha,r,\eta)$ independent of $h$ such that
\begin{equation}
\Vert y-y_h\Vert_{L^2(Q_T)} \leq  K \frac{h^k}{\sqrt{\eta}}. 
\end{equation}
\end{theorem}

\section{Numerical experiments}\label{sec_experiment}

We now report and discuss some numerical experiments corresponding to mixed formulation (\ref{eq:mfh}) and (\ref{eq:mfalphah}) for $N=1$ and $N=2$.

\subsection{One dimensional case ($N = 1$)}

We take $\Omega=(0,1)$. In order to check the convergence of the method, we consider explicit solutions of (\ref{eq:wave}). We define the smooth initial condition (see \cite{cindea_moireau}): 
\begin{equation}
\nonumber
(\textbf{EX1})\quad
\left\{
\begin{aligned}
& y_0(x)=16x^2(1-x)^2, \\
& y_1(x)= (3x-4x^3)\, 1_{(0,0.5)}(x) + (4x^3-12x^2+9x-1)\, 1_{(0.5,1)}(x), 
\end{aligned}
x\in (0,1)
\right.
\end{equation}
and $f=0$. The corresponding solution of (\ref{eq:wave}) with $c\equiv 1, d\equiv 0$ is given by 
$$
y(x,t)=\sum_{k>0}  \biggl(a_k \cos(k\pi  t) + \frac{b_k}{k\pi}\sin(k\pi t)\biggr) \sqrt{2} \sin(k\pi x)
$$ 
with 
$$
a_k=\frac{32\sqrt{2}(\pi^2k^2-12)}{\pi^5 k^5} ((-1)^k-1), \quad b_k=\frac{48\sqrt{2}\sin(\pi k/2)}{\pi^4 k^4}, \quad k>0.
$$
We also define the initial data in $H_0^1(\Omega)\times L^2(\Omega)$
\begin{equation}
\nonumber
(\textbf{EX2}) \quad
y_0(x)=1-\vert 2x-1\vert, \quad y_1(x)=  1_{(1/3,2/3)}(x),  \qquad x\in (0,1)
\end{equation}
for which the Fourier coefficients are 
$$
a_k=\frac{4\sqrt{2}}{\pi^2 k^2} \sin(\pi k/2), \quad b_k=  \frac{1}{\pi k}(\cos(\pi k/3)-\cos(2\pi k/3)), \quad k>0.
$$

\subsubsection{The cylindrical case: $q_T=\omega\times (0,T)$}
\label{sec:eps0}

We consider the case $\eps=0$  described in Section \ref{recovering_y}. We take $\omega=(0.1,0.3)$ and $T=2$ for which the inequality (\ref{iobs}) holds true. 
We consider the BFS finite element with uniform triangulation (each element $K$ of the triangulation $\mathcal{T}_h$ is a rectangle of lengths $\Delta x$ and $\Delta t$ so that $h=\sqrt{(\Delta x)^2+(\Delta t)^2}$). We recall that the direct method amounts to solve, for any $h$, the linear system (\ref{matrixmfh}). We use the LU decomposition method. Table \ref{tab:ex1_T2} collects some norms with respect to $h$ for the initial data (\textbf{EX1}) for $r=1$ and for $\Delta x=\Delta t$. We observe a linear convergence for the variables $y_h$, $\lambda_h$ for the $L^2$-norm:  

\begin{equation}
\frac{\Vert y-y_h\Vert_{L^2(Q_T)}}{\Vert y\Vert_{L^2(Q_T)}}=\mathcal{O}(h^{1.03}), \quad \frac{\Vert y-y_h\Vert_{L^2(q_T)}}{\Vert y\Vert_{L^2(q_T)}}=\mathcal{O}(h^{0.98}), \quad \Vert \lambda_h\Vert_{L^2(Q_T)}=\mathcal{O}(h^{0.98}).
\end{equation}
In agreement with Remark \ref{rk_lambda_sys}, since $y_{obs}$ is by construction the restriction to $q_T$ of a solution of (\ref{eq:wave}), the sequence $\lambda_h$, approximation of $\lambda$, vanishes as $h\to 0$. 
The $L^2$-norm of $Ly_h$ do also converges to $0$ with $h$, with a lower rate: 
\begin{equation}
\Vert L y_h\Vert_{L^2(Q_T)}=\mathcal{O}(h^{0.71}).
\end{equation}

\begin{table}[http]
\centering
\begin{tabular}{|c|ccccc|}
\hline
$h$  &  $7.01\times 10^{-2}$  & $3.53\times 10^{-2}$ & $1.76\times 10^{-2}$ & $8.83\times 10^{-3}$ & $4.42\times 10^{-3}$  \tabularnewline

\hline 
$\frac{\Vert y-y_h\Vert_{L^2(Q_T)}}{\Vert y\Vert_{L^2(Q_T)}}$ & $9.55\times 10^{-2}$ & $4.58\times 10^{-2}$ & $2.24\times 10^{-2}$ & $1.10\times 10^{-2}$ & $5.52\times 10^{-3}$\tabularnewline

$\frac{\Vert y-y_h\Vert_{L^2(q_T)}}{\Vert y\Vert_{L^2(q_T)}}$ & $8.35\times 10^{-2}$ & $4.28\times 10^{-2}$ & $2.16\times 10^{-2}$ & $1.09\times 10^{-2}$ & $5.51\times 10^{-3}$ \tabularnewline

$\Vert L y_h\Vert_{L^2(Q_T)}$ & $5.62\times 10^{-3}$ & $3.21\times 10^{-3}$ & $1.78\times 10^{-3}$ & $9.99\times 10^{-4}$ & $8.54\times 10^{-4}$\tabularnewline

$\Vert \lambda_h\Vert_{L^2(Q_T)}$ & $2.67\times 10^{-5}$ & $1.37\times 10^{-5}$ & $6.89\times 10^{-6}$ & $3.44\times 10^{-6}$  & $1.76\times 10^{-6}$\tabularnewline

$\kappa$ & $1.4\times 10^{10}$ & $4.6\times 10^{11}$ & $1.3\times 10^{13}$ & $4.2\times 10^{14}$ &  $1.3\times 10^{16}$ \tabularnewline

card($\{\lambda_h\}$) & $861$ & $3\ 321$ & $13\ 041$ & $51\ 681$ & $205\ 761$\tabularnewline

$\sharp$ \textrm{CG iterates} & $27$ & $42$ & $70$ & $96$ & $90$\tabularnewline

\hline
\end{tabular}
\caption{Example \textbf{EX1} - $r=1$ - $T=2$ - $\Vert y\Vert_{L^2(q_T)}=5.95\times 10^{-2}$ -  $\Vert y\Vert_{L^2(Q_T)}=1.59\times 10^{-1}$.}
\label{tab:ex1_T2}
\end{table}

We also check that the minimization of the functional $J^{\star\star}_r$ introduced in Proposition \ref{prop_equiv_dual} leads exactly to the same result: we recall that the minimization of the functional $J^{\star\star}_r$ corresponds to the resolution of the associate mixed formulation by an iterative Uzawa type method. The minimization is done using a conjugate gradient algorithm ( we refer to \cite[Section 2.2]{NC-AM-mixedwave} for the algorithm). Each iteration amounts to solve a linear system involving the matrix $A_{r,h}$ which is sparse, symmetric and positive definite. The Cholesky method is used. The performance of the algorithm depends on the conditioning number of the operator $\mathcal{P}_r$: precisely, it is known that (see for instance \cite{Daniel1971}), 
\[
\Vert \lambda^n-\lambda\Vert_{L^2(Q_T)} \leq 2\sqrt{\nu(\mathcal{P}_r)}   \biggl(\frac{\sqrt{\nu(\mathcal{P}_r)}-1}{\sqrt{\nu(\mathcal{P}_r)}+1}  \biggr)^n \Vert \lambda^0-\lambda\Vert_{L^2(Q_T)}, \quad \forall n\geq 1
\]
where $\lambda$ minimizes $J_r^{\star\star}$. $\nu(\mathcal{P}_r)=\Vert \mathcal{P}_r\Vert \Vert \mathcal{P}_r^{-1}\Vert$ denotes the condition number of the operator $\mathcal{P}_r$. As discussed in \cite[Section 4.4]{NC-AM-mixedwave}, the conditioning number of $\mathcal{P}_r$ restricted to $\Lambda_h\subset L^2(Q_T)$ behaves asymptotically as $C_r^{-2}h^{-2}$. Table \ref{tab:ex1_T2} reports the number of iterations of the algorithm, initiated with $\lambda^0=0$ in $Q_T$. We take $\epsilon=10^{-10}$ as a stopping threshold for the algorithm (the algorithm is stopped as soon as the norm of the residue $g^n$ given here by $Ly^n$ satisfies $\Vert g^n\Vert_{L^2(Q_T)}\leq 10^{-10} \Vert g^0\Vert_{L^2(Q_T)}$).

Table \ref{tab:ex1_T2} reports the number of iterates to reach convergence, with respect to $h$. We observe that this number is sub-linear with respect to $h$, precisely, with respect to the dimension $m_h=card(\{\lambda_h\})$ of the approximated problems. This renders this method very attractive from a numerical point of view. 

From Remark \ref{FVeps}, we also check the convergence w.r.t. $h$ when we assume from the beginning that the multiplier $\lambda$ vanishes  (see Table \ref{tab:ex1_T2_lambda}). This amounts to minimize the functional $J_r$
given by  (\ref{def_Jyr}) or, equivalently, to perform exactly one iteration of the conjugate gradient algorithm we have just discussed. With $r=1$, we observe a weaker convergence : 
\begin{equation}
\frac{\Vert y-y_h\Vert_{L^2(Q_T)}}{\Vert y\Vert_{L^2(Q_T)}}=\mathcal{O}(h^{0.574}), \quad \frac{\Vert y-y_h\Vert_{L^2(q_T)}}{\Vert y\Vert_{L^2(q_T)}}=\mathcal{O}(h^{0.94}).
\end{equation}
This example illustrates that the convergence of  $Ly_h$ to $0$ in the norm $L^2(0,T,H^{-1}(0,1))$ is enough here to guarantee the convergence of the approximation $y_h$: we get 
that  $h \Vert L y_h\Vert_{L^2(Q_T)} \approx \Vert L y_h\Vert_{L^2(0,T; H^{-1}(0,1)}=\mathcal{O}(h^{0.3})$ while $\Vert L y_h\Vert_{L^2(Q_T)}$ slightly increases. Obviously, in this specific situation, a larger $r$ (acting as a penalty term) independent of $h$ yields a lower $\Vert Ly_h\Vert_{L^2(Q_T)}$ norm.  

\begin{table}[http]
\centering
\begin{tabular}{|c|ccccc|}
\hline
$h$  &  $7.01\times 10^{-2}$  & $3.53\times 10^{-2}$ & $1.76\times 10^{-2}$ & $8.83\times 10^{-3}$ & $4.42\times 10^{-3}$  \tabularnewline

\hline $\frac{\Vert y-y_h\Vert_{L^2(Q_T)}}{\Vert y\Vert_{L^2(Q_T)}}$ & $9.74\times 10^{-2}$ & $4.90\times 10^{-2}$ & $2.84\times 10^{-2}$ & $2.16\times 10^{-2}$ & $2.01\times 10^{-2}$\tabularnewline

$\frac{\Vert y-y_h\Vert_{L^2(q_T)}}{\Vert y\Vert_{L^2(q_T)}}$ & $8.35\times 10^{-2}$ & $4.28\times 10^{-2}$ & $2.18\times 10^{-2}$ & $1.12\times 10^{-2}$ & $6.21\times 10^{-3}$\tabularnewline

$\Vert L y_h\Vert_{L^2(Q_T)}$ & $7.72\times 10^{-3}$ & $1.11\times 10^{-2}$ & $2.01\times 10^{-2}$ & $3.40\times 10^{-2}$ & $4.79\times 10^{-2}$\tabularnewline


\hline
\end{tabular}
\caption{Example \textbf{EX1} - $r=1$ - $T=2$ - $\lambda$ fixed to zero.}
\label{tab:ex1_T2_lambda}
\end{table}

On the contrary, we check that the convergence to $0$ of $\Vert y-y_h\Vert_{L^2(Q_T)}$ is lost when the inequality (\ref{iobs}) is not satisfied: Table \ref{tab:ex1_T1} collects the norms w.r.t. $h$ for the same data except the value $T=1$ (for which the uniqueness of the solution is lost): we observe that  $\Vert y-y_h\Vert_{L^2(Q_T)}$ increases as $h\to 0$.
As an illustration of the loss of uniqueness, these value also yields to a larger conditioning number $\kappa$ of the matrix $A_{r,h}$.


\begin{table}[http]
\centering
\begin{tabular}{|c|ccccc|}
\hline
$h$  &  $7.01\times 10^{-2}$  & $3.53\times 10^{-2}$ & $1.76\times 10^{-2}$ & $8.83\times 10^{-3}$ & $4.42\times 10^{-3}$  \tabularnewline

\hline 
$\frac{\Vert y-y_h\Vert_{L^2(Q_T)}}{\Vert y\Vert_{L^2(Q_T)}}$ & $1.21\times 10^{-1}$ & $1.08\times 10^{-1}$ & $1.34\times 10^{-1}$ & $2.42\times 10^{-1}$ & $5.19\times 10^{-1}$\tabularnewline

$\frac{\Vert y-y_h\Vert_{L^2(q_T)}}{\Vert y\Vert_{L^2(q_T)}}$ & $8.40\times 10^{-2}$ & $4.34\times 10^{-2}$ & $2.22\times 10^{-2}$ & $1.12\times 10^{-2}$ & $5.62\times 10^{-3}$\tabularnewline

$\Vert L y_h\Vert_{L^2(Q_T)}$ & $5.62\times 10^{-2}$ & $2.77\times 10^{-2}$ & $2.63\times 10^{-2}$ & $2.25\times 10^{-2}$ & $2.15\times 10^{-2}$\tabularnewline

$\Vert \lambda_h\Vert_{L^2(Q_T)}$ & $1.84\times 10^{-5}$ & $9.48\times 10^{-6}$ & $4.76\times 10^{-6}$ & $2.38\times 10^{-6}$  & $1.19\times 10^{-6}$\tabularnewline

$\kappa$ & $1.2\times 10^{11}$ & $9.8\times 10^{12}$ & $1.1\times 10^{15}$ & $1.5\times 10^{17}$ &  $2.7\times 10^{19}$ \tabularnewline

\hline
\end{tabular}
\caption{Example \textbf{EX1} - $r=1$ - $T=1$ - $\Vert y_{ex}\Vert_{L^2(q_T)}=4.21\times 10^{-2}$ -  $\Vert y_{ex}\Vert_{L^2(Q_T)}=1.12\times 10^{-1}$.}
\label{tab:ex1_T1}
\end{table}

Similar conclusions hold with the less regular initial data (\textbf{EX2}). Numerical results are reported in Table \ref{tab:ex2_T2}. We still observe a linear convergence w.r.t. $h$ of 
$\Vert y-y_h\Vert_{L^2(Q_T)}$, $\Vert y-y_h\Vert_{L^2(q_T)}$ and $\Vert \lambda_h\Vert_{L^2(Q_T)}$. One notable difference is that the convergence rate is weaker for the norm $\Vert L y_h\Vert_{L^2(Q_T)}$: 
\begin{equation}
\Vert L y_h\Vert_{L^2(Q_T)}=\mathcal{O}(h^{0.123}).
\end{equation}
Again, this is enough to guarantee the convergence of $y_h$ toward a solution of the wave equation: recall that then $\Vert L y_h\Vert_{L^2(0,T;H^{-1}(0,1))}=\mathcal{O}(h^{1.123})$.
We also observe that the number of iterates in the $CG$ algorithm remains largely sub-linear but is slightly larger: precisely, we have $\sharp$ iter $=\mathcal{O}(h^{-0.71})$. 
Table \ref{tab:ex2_T1} illustrates the case $T=1$ while Table \ref{tab:ex2_T2_lambda} illustrates the minimization of $J_r$ (see \ref{def_Jyr}), both for $r=1$.

\begin{table}[http]
\centering
\begin{tabular}{|c|ccccc|}
\hline
$h$  &  $7.01\times 10^{-2}$  & $3.53\times 10^{-2}$ & $1.76\times 10^{-2}$ & $8.83\times 10^{-3}$ & $4.42\times 10^{-3}$  \tabularnewline

\hline

$\frac{\Vert y-y_h\Vert_{L^2(Q_T)}}{\Vert y\Vert_{L^2(Q_T)}}$ & $1.01\times 10^{-1}$ & $4.81\times 10^{-2}$ & $2.34\times 10^{-2}$ & $1.15\times 10^{-2}$ & $5.68\times 10^{-3}$\tabularnewline

$\frac{\Vert y-y_h\Vert_{L^2(q_T)}}{\Vert y\Vert_{L^2(q_T)}}$ & $1.34\times 10^{-1}$ & $5.05\times 10^{-2}$ & $2.37\times 10^{-2}$ & $1.16\times 10^{-2}$ & $5.80\times 10^{-3}$\tabularnewline

$\Vert L y_h\Vert_{L^2(Q_T)}$ & $7.18\times 10^{-2}$ & $6.59\times 10^{-2}$ & $6.11\times 10^{-2}$ & $5.55\times 10^{-2}$ & $5.10\times 10^{-2}$\tabularnewline

$\Vert \lambda_h\Vert_{L^2(Q_T)}$ & $1.07\times 10^{-4}$ & $4.70\times 10^{-5}$ & $2.32\times 10^{-5}$ & $1.15\times 10^{-5}$  & $5.76\times 10^{-6}$\tabularnewline

$\sharp$ \textrm{CG iterates} & $29$ & $46$ & $83$ & $133$ & $201$\tabularnewline
    
\hline
\end{tabular}
\caption{Example \textbf{EX2} - $r=1$ - $T=2$ - $\Vert y\Vert_{L^2(q_T)}=1.56\times 10^{-1}$ -  $\Vert y\Vert_{L^2(Q_T)}=4.14\times 10^{-1}$.}
\label{tab:ex2_T2}
\end{table}

\begin{table}[http]
\centering
\begin{tabular}{|c|ccccc|}
\hline
$h$  &  $7.01\times 10^{-2}$  & $3.53\times 10^{-2}$ & $1.76\times 10^{-2}$ & $8.83\times 10^{-3}$ & $4.42\times 10^{-3}$  \tabularnewline

\hline 

$\frac{\Vert y-y_h\Vert_{L^2(Q_T)}}{\Vert y\Vert_{L^2(Q_T)}}$ & $2.74\times 10^{-1}$ & $4.15\times 10^{-1}$ & $6.30\times 10^{-1}$ & $1.21$ & $2.62$\tabularnewline

$\frac{\Vert y-y_h\Vert_{L^2(q_T)}}{\Vert y\Vert_{L^2(q_T)}}$ & $1.37\times 10^{-1}$ & $5.76\times 10^{-2}$ & $2.89\times 10^{-2}$ & $2.41\times 10^{-2}$ & $7.76\times 10^{-3}$\tabularnewline

$\Vert L y_h\Vert_{L^2(Q_T)}$ & $5.97\times 10^{-2}$ & $4.96\times 10^{-2}$ & $4.96\times 10^{-2}$ & $4.52\times 10^{-2}$ & $4.21\times 10^{-2}$\tabularnewline

$\Vert \lambda_h\Vert_{L^2(Q_T)}$ & $4.97\times 10^{-5}$ & $2.32\times 10^{-5}$ & $1.15\times 10^{-5}$ & $5.76\times 10^{-5}$  & $2.87\times 10^{-6}$\tabularnewline

\hline
\end{tabular}
\caption{Example \textbf{EX2} - $r=1$ - $T=1$ - $\Vert y\Vert_{L^2(q_T)}=1.104\times 10^{-1}$ -  $\Vert y\Vert_{L^2(Q_T)}=2.93\times 10^{-1}$.}
\label{tab:ex2_T1}
\end{table}

\begin{table}[http]
\centering
\begin{tabular}{|c|ccccc|}
\hline
$h$  &  $7.01\times 10^{-2}$  & $3.53\times 10^{-2}$ & $1.76\times 10^{-2}$ & $8.83\times 10^{-3}$ & $4.42\times 10^{-3}$  \tabularnewline

\hline 

$\frac{\Vert y-y_h\Vert_{L^2(Q_T)}}{\Vert y\Vert_{L^2(Q_T)}}$ & $1.02\times 10^{-1}$ & $5.27\times 10^{-2}$ & $3.18\times 10^{-2}$ & $2.48\times 10^{-2}$ & $2.25\times 10^{-2}$\tabularnewline

$\frac{\Vert y-y_h\Vert_{L^2(q_T)}}{\Vert y\Vert_{L^2(q_T)}}$ & $1.34\times 10^{-1}$ & $5.06\times 10^{-2}$ & $2.37\times 10^{-2}$ & $1.21\times 10^{-2}$ & $6.65\times 10^{-3}$\tabularnewline

$\Vert L y_h\Vert_{L^2(Q_T)}$ & $7.43\times 10^{-2}$ & $7.43\times 10^{-2}$ & $8.65\times 10^{-2}$ & $1.10\times 10^{-1}$ & $1.37\times 10^{-2}$\tabularnewline

\hline
\end{tabular}
\caption{Example \textbf{EX2} - $r=1$ - $T=2$ - $\lambda$ fixed to zero.}
\label{tab:ex2_T2_lambda}
\end{table}

We end this section with some numerical results for the stabilized mixed formulation (\ref{eq:mfalphah}). The main difference is that the multiplier $\lambda$
is approximated in a much richer space $\widetilde\Lambda_h$ (see \ref{choice_lambdah_alpha}) leading to larger linear system. Table \ref{tab:ex2_T2_stab}
consider the case of the example \textbf{EX2} with $T=2$ and $\alpha=1/2$. In order to compare with the formulation (\ref{eq:mfh}), we take again $r=1$.  We observe the convergence w.r.t $h$ and obtain slightly better rates and constants than in Table \ref{tab:ex2_T2}: in particular, we have $\Vert y-y_h\Vert_{L^2(Q_T)} / \Vert y\Vert_{L^2(Q_T)} =\mathcal{O}(h^{1.10})$. 
This is partially due to the fact that the space $\widetilde{\Lambda}_h$ used for the variable $\lambda_h$ in (\ref{eq:mfalphah}) is richer than the space $\Lambda_h$ used in (\ref{eq:mfh}). However, for $\alpha=0$ leading to the non stabilized mixed formulation, the space $\widetilde{\Lambda}_h$
is too rich and produce poor result, while we obtain very similar results for any values of $\alpha$ in $(0,1]$. Finally, we also check that - in contrast with the mixed formulation (\ref{eq:mfh}) - the positive parameter $r$ does not affect the numerical results. 

\begin{table}[http]
\centering
\begin{tabular}{|c|ccccc|}
\hline
$h$  &  $7.01\times 10^{-2}$  & $3.53\times 10^{-2}$ & $1.76\times 10^{-2}$ & $8.83\times 10^{-3}$ & $4.42\times 10^{-3}$  \tabularnewline

\hline 
$\frac{\Vert y-y_h\Vert_{L^2(Q_T)}}{\Vert y\Vert_{L^2(Q_T)}}$ & $8.48\times 10^{-2}$ & $4.01\times 10^{-2}$ & $1.85\times 10^{-2}$ & $8.66\times 10^{-3}$ & $4.01\times 10^{-3}$\tabularnewline

$\frac{\Vert y-y_h\Vert_{L^2(q_T)}}{\Vert y\Vert_{L^2(q_T)}}$ & $2.80\times 10^{-1}$ & $7.26\times 10^{-2}$ & $2.61\times 10^{-2}$ & $1.12\times 10^{-2}$ & $5.05\times 10^{-3}$ \tabularnewline

$\Vert L y_h\Vert_{L^2(Q_T)}$ & $7.25\times 10^{-2}$ & $6.59\times 10^{-2}$ & $6.16\times 10^{-2}$ & $5.58\times 10^{-2}$ & $5.08\times 10^{-2}$\tabularnewline

$\Vert \lambda_h\Vert_{L^2(Q_T)}$ & $4.11\times 10^{-3}$ & $2.04\times 10^{-3}$ & $1.49\times 10^{-3}$ & $1.01\times 10^{-3}$  & $7.37\times 10^{-4}$\tabularnewline

\hline
\end{tabular}
\caption{Example \textbf{EX2} - $r=1$ - $T=2$ - $\alpha=1/2$ - $\Vert y\Vert_{L^2(q_T)}=5.95\times 10^{-2}$ -  $\Vert y\Vert_{L^2(Q_T)}=1.59\times 10^{-1}$.}
\label{tab:ex2_T2_stab}
\end{table}

\begin{figure}[ht!]
\centering
\begin{tabular}{cc}
\includegraphics[width=0.3\textwidth]{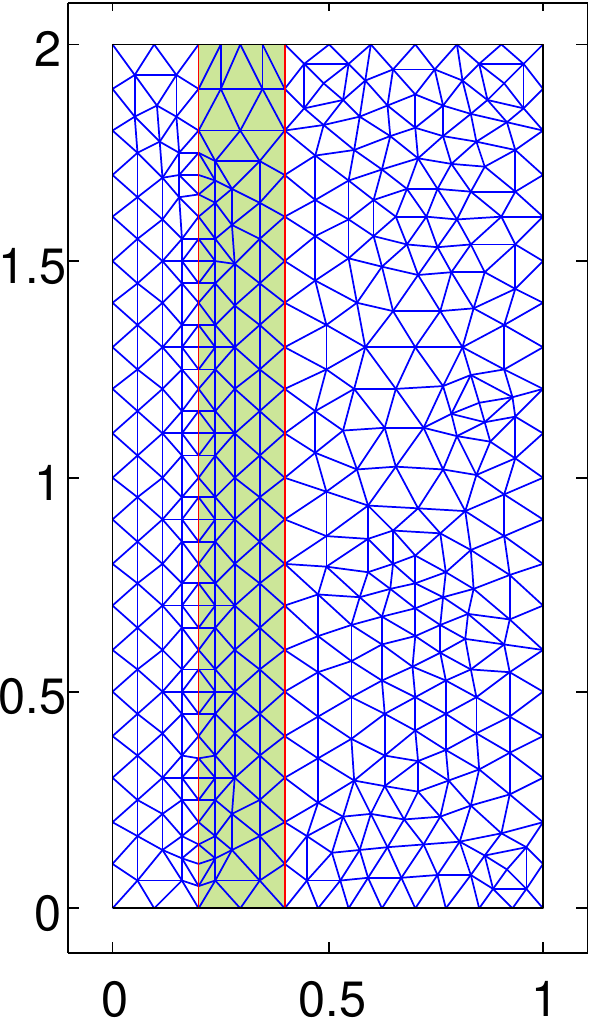} & \hspace{1cm}\includegraphics[width=0.3\textwidth]{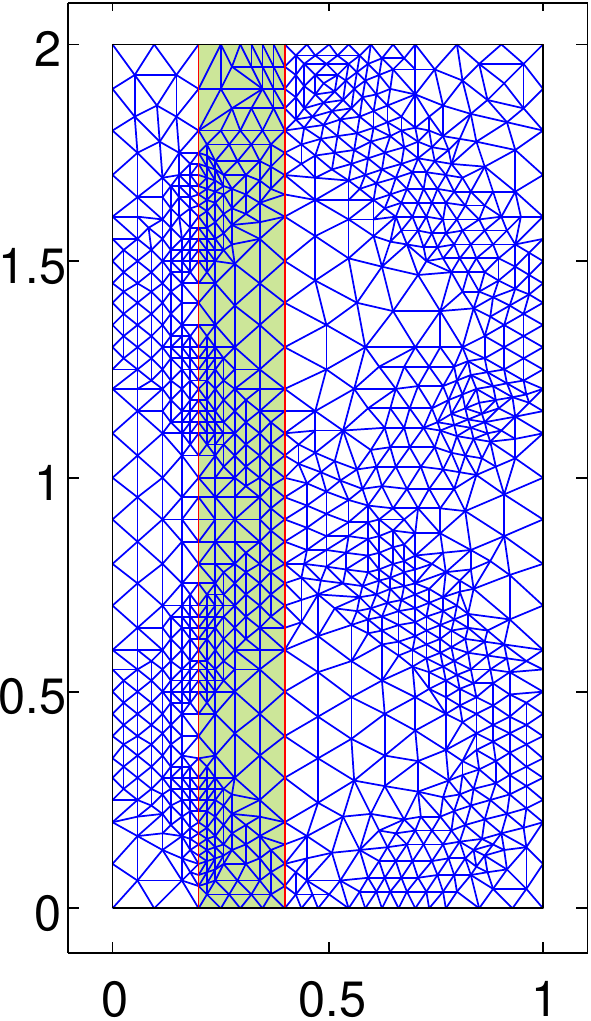} \\
\includegraphics[width=0.3\textwidth]{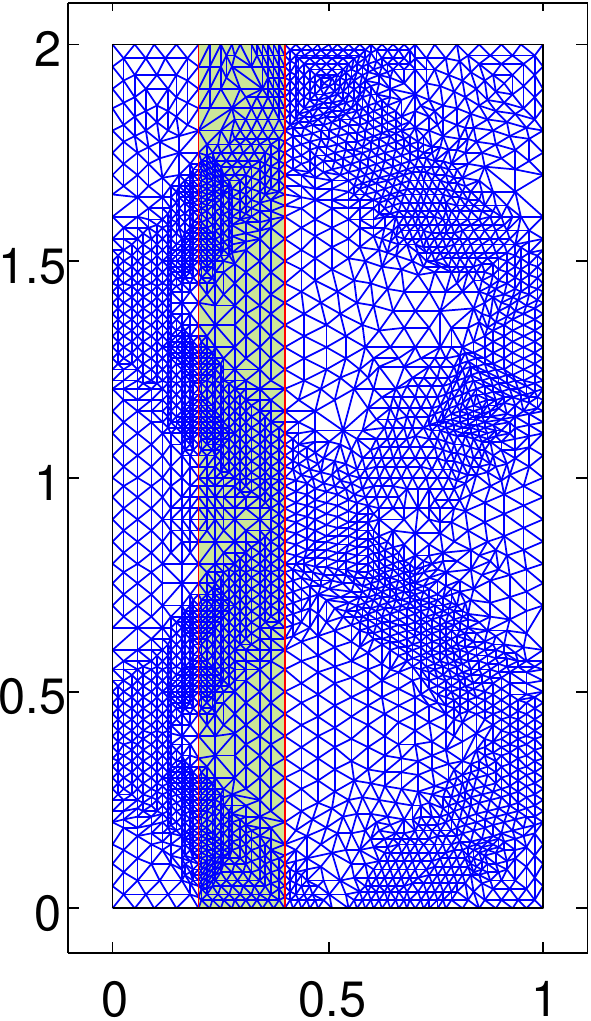} & \hspace{1cm}\includegraphics[width=0.3\textwidth]{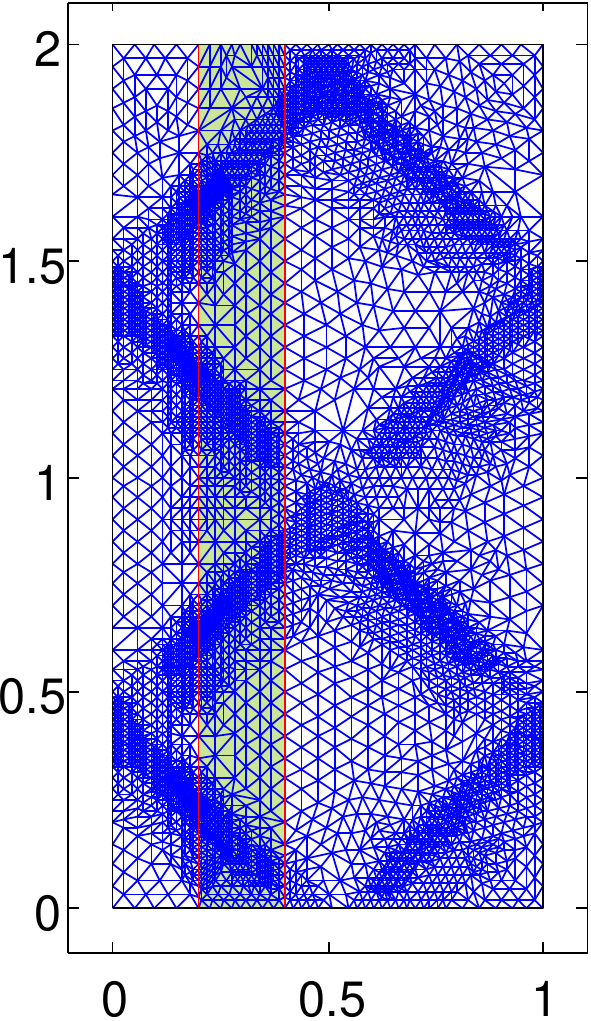} 
\end{tabular}
\caption{Iterative refinement of the triangular mesh over $Q_T$ with respect to the variable $y$.}
\label{fig:adapt}
\end{figure}

We also emphasize that this variational method which requires a finite element discretization of the time-space $Q_T$ is particularly well-adapted to mesh optimization. Still for the example \textbf{EX2}, Figure \ref{fig:adapt} depicts a sequence of four distinct meshes of $Q_T=(0,1)\times (0,T)$: the sequence is initiated with a coarse and regularly distributed mesh. The three other meshes are successively obtained by local refinement based on the norm of the gradient of $y_h$ on each triangle of $\mathcal{T}_h$. As expected, the refinement is concentrated around the lines of singularity of $y_h$ travelling in $Q_T$, generated by the singularity of the initial position $y_0$. The four meshes contain  $792, 2\ 108, 7\ 902$ and $14\ 717$ triangles respectively (see Table \ref{tab:adapt}). The results obtained using the reduced HCT finite element are reported in Table \ref{tab:adapt}.

\begin{table}[ht]
\centering
\begin{tabular}{|c|cccc|}
\hline
Mesh number  &  1 & 2 & 3 & 4  \tabularnewline
\hline 
$\sharp$ elements & $792$ & $2\ 108$ & $7\ 902$ & $14\ 717$ \tabularnewline
$\sharp$ points & $429$ & $1\ 101$ & $4\ 041$ & $7\ 462$ \tabularnewline
$\frac{\Vert y-y_h\Vert_{L^2(Q_T)}}{\Vert y\Vert_{L^2(Q_T)}}$ & $1.34 \times 10^{-2}$ & $8.69 \times 10^{-3}$ & $6.01 \times 10^{-3}$ & $5.9 \times 10^{-3}$ \tabularnewline
$\Vert\lambda \Vert_{L^2(Q_T)}$ & $1.14 \times 10^{-5}$ & $7.99 \times 10^{-6}$ & $5.02 \times 10^{-6}$ & $4.79 \times 10^{-6}$\tabularnewline
\hline
\end{tabular}
\caption{Example \textbf{(EX2)} - Information concerning the meshes and approximation errors for mesh adaptation strategy.}
\label{tab:adapt}
\end{table}

\subsubsection{The non-cylindrical case}

\begin{figure}[ht]
\centering
\begin{tabular}{ccc}
\includegraphics[width=0.28\textwidth]{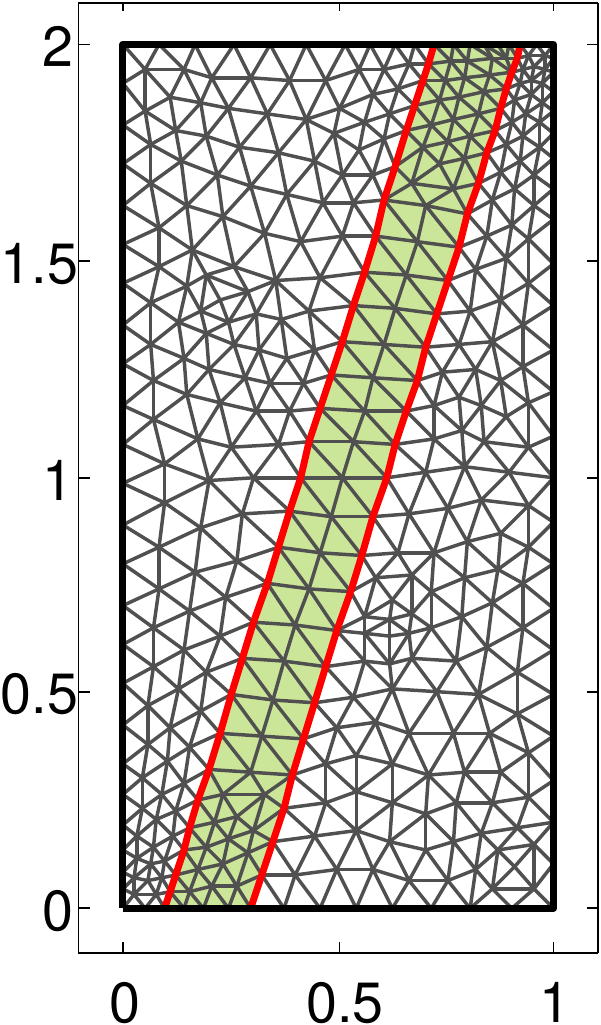} & \qquad & 
\includegraphics[width=0.28\textwidth]{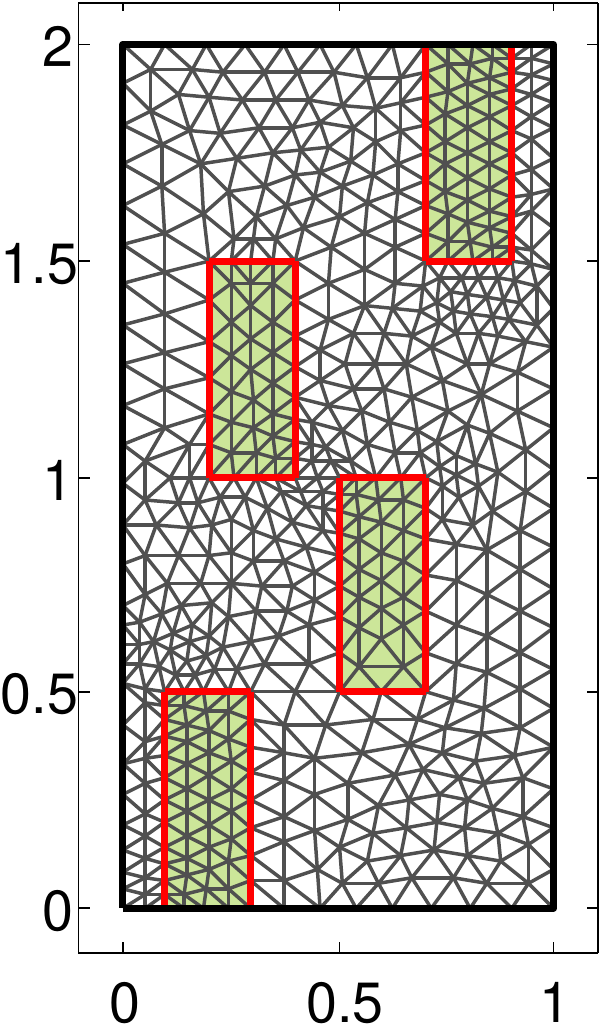} \\
(a) & &(b)
\end{tabular}
\caption{Domain $q_T^1$ (a) and domain $q_T^2$ (b) triangulated using some coarse meshes.}
\label{fig:qTvar}
\end{figure}

We numerically illustrate the reconstruction of the state of the wave equation \eqref{eq:wave} from measurements $y_{obs}$ which are available in domains $q_T \subset Q_T$ non-constant in time (considered recently in \cite{CC-NC-AM} in a controllability context). Time dependent domains also appears for time under sampled observations (or measurements): we refer to \cite{cindea_imperiale_moireau}.  In what follows we take $T = 2$ and $q_T$ to be  one of the two following domains:
\begin{equation}\label{eq:qTvar1}
q_{T}^1 := \left\{ 
(x, t) \in Q_T \text{ such that } \left| x - \frac{3t}{5T} - \frac{1}{5}\right| < \frac{1}{10} \text{ for every } t \in (0, T)
\right\},
\end{equation}
\begin{align} \nonumber
q_T^2  := & \left( \frac{1}{10}, \frac{2}{10}\right) \times \left(0, \frac{T}{4}\right) \bigcup \left( \frac{1}{2}, \frac{7}{10}\right) \times \left(\frac{T}{4}, \frac{T}{2}\right) \\
\label{eq:qTvar2}
& \bigcup \left( \frac{1}{5}, \frac{2}{5}\right) \times \left(\frac{T}{2}, \frac{3T}{4}\right) \bigcup \left( \frac{7}{10}, \frac{9}{10}\right) \times \left(\frac{3T}{4}, T\right).
\end{align}
These two pairs $(T,q^i_T)$ $i=1,2$ satisfy the standard geometric optic condition: therefore, using \cite{CC-NC-AM}, Proposition 2.1, inequality (\ref{iobs}) holds true. Both domains $q_T^1$ and $q_T^2$ are displayed in Figure \ref{fig:qTvar} with the coarsest of the meshes that are used for the numerical experiments in this section. 

We consider five levels of regular triangular meshes and use the reduced Hsieh-Clough-Tocher finite element.
We illustrate our method on the reconstruction of the solution of the wave equation corresponding to initial data \textbf{(EX2)} considered in Section \ref{sec:eps0}.

Since domains $q_T^1$ and $q_T^2$ satisfy the geometric optic condition, we obtain similar results as in the case $q_T = \omega \times (0, T)$ studied in the previous section. More precisely, these results are reported in Table \ref{tab:ex2_T2_qT1} and Table \ref{tab:ex2_T2_qT2} for domain $q_T^1$ and $q_T^2$ respectively.

\begin{table}[http]
\centering
\begin{tabular}{|c|ccccc|}
\hline
$h$  &   $7.18 \times 10^{-2}$ & $3.59 \times 10^{-2}$ & $1.79 \times 10^{-2}$ & $9 \times 10^{-3}$ & $4.5 \times 10^{-3}$ \tabularnewline

\hline

$\frac{\Vert y-y_h\Vert_{L^2(Q_T)}}{\Vert y\Vert_{L^2(Q_T)}}$ & $ 2.02\times 10^{-2} $ & $ 7.83\times 10^{-3} $ & $ 3.32\times 10^{-3} $ & $ 1.36\times 10^{-3} $ & $ 6.27\times 10^{-4} $\tabularnewline

$\frac{\Vert y-y_h\Vert_{L^2(q_T)}}{\Vert y\Vert_{L^2(q_T)}}$ & $ 1.85\times 10^{-2} $ & $ 6.69\times 10^{-3} $ & $ 2.40\times 10^{-3} $ & $ 1.03\times 10^{-3} $ & $ 4.56\times 10^{-4} $\tabularnewline


$\Vert L y_h\Vert_{L^2(Q_T)}$ & $ 3.41 $ & $ 3.78 $ & $ 4.15 $ & $ 4.47 $ & $ 4.76 $ \tabularnewline

$\Vert \lambda_h\Vert_{L^2(Q_T)}$ & $ 1.97\times 10^{-5} $ & $ 7.03\times 10^{-6} $ & $ 1.70\times 10^{-6} $ & $ 4.14\times 10^{-7} $ & $ 1.10\times 10^{-7} $\tabularnewline

$\kappa$ & $ 1.18\times 10^{8} $ & $ 1.84\times 10^{9} $ & $ 1.61\times 10^{10} $ & $ 1.75\times 10^{11} $ & $ 1.38\times 10^{12} $ \tabularnewline

$\text{card}(\{ \lambda_h \})$ & $ 429 $ & $ 1\ 633 $ & $ 6\ 369 $ & $ 25\ 153 $ & $ 99\ 969 $ \tabularnewline

$\sharp$ \textrm{CG iterates} & $108$ & $206$ & $392$ & $954$ & $2\ 009$\tabularnewline
    
\hline
\end{tabular}
\caption{Observation domain $q_T^1$. Example \textbf{EX2} - $r=1$ - $T=2$ - $\Vert y\Vert_{L^2(q_T)}=2.75\times 10^{-1}$ -  $\Vert y\Vert_{L^2(Q_T)}=5.87\times 10^{-1}$.}
\label{tab:ex2_T2_qT1}
\end{table}

Remark that the number of iterations needed for the conjugate gradient algorithm in order to achieve a residual smaller than $10^{-10}$ when we minimize the functional $J^{\star\star}$ over $\Lambda_h$ is slightly larger than in the situations described in the previous section. 

\begin{table}[http]
\centering
\begin{tabular}{|c|ccccc|}
\hline
$h$  &   $6.24 \times 10^{-2}$ & $3.12 \times 10^{-2}$ & $1.56 \times 10^{-2}$ & $7.8 \times 10^{-3}$ & $3.9 \times 10^{-3}$ \tabularnewline

\hline

$\frac{\Vert y-y_h\Vert_{L^2(Q_T)}}{\Vert y\Vert_{L^2(Q_T)}}$ & $ 1.38\times 10^{-2} $ & $ 6.37\times 10^{-3} $ & $ 2.64\times 10^{-3} $ & $ 1.15\times 10^{-3} $ & $ 5.25\times 10^{-4} $   \tabularnewline

$\frac{\Vert y-y_h\Vert_{L^2(q_T)}}{\Vert y\Vert_{L^2(q_T)}}$ & $ 1.27\times 10^{-2} $ & $ 4.79\times 10^{-3} $ & $ 2.02\times 10^{-3} $ & $ 9.11\times 10^{-4} $ & $ 4.29\times 10^{-4} $ \tabularnewline

$\Vert L y_h\Vert_{L^2(Q_T)}$ & $ 3.86 $ & $ 3.45 $ & $ 3.36 $ & $ 3.85 $ & $ 4.16 $ \tabularnewline

$\Vert \lambda_h\Vert_{L^2(Q_T)}$ & $ 6.37\times 10^{-6} $ & $ 1.65 \times 10^{-6} $ & $ 3.88\times 10^{-7} $ & $ 9.74\times 10^{-8} $ & $ 2.90\times 10^{-8} $ \tabularnewline

$\kappa$ & $ 2.02\times 10^{8} $ & $ 2.62\times 10^{9} $ & $ 2.05\times 10^{10} $ & $ 1.61\times 10^{11} $ & $ 1.32\times 10^{12} $ \tabularnewline

$\text{card}(\{ \lambda_h \})$ & $ 554 $ & $ 2\ 135 $ & $ 8\ 381 $ & $ 33\ 209 $ & $ 132\ 209 $\tabularnewline

$\sharp$ \textrm{CG iterates} & $141$ & $331$ & $720$ & $1\ 446$ & $3\ 318$\tabularnewline
    
\hline
\end{tabular}
\caption{Observation domain $q_T^2$. Example \textbf{EX2} - $r=1$ - $T=2$ - $\Vert y\Vert_{L^2(q_T)}=2.75\times 10^{-1}$ -  $\Vert y\Vert_{L^2(Q_T)}=5.87\times 10^{-1}$.}
\label{tab:ex2_T2_qT2}
\end{table}

\begin{figure}[ht]
\begin{tabular}{cc}
\includegraphics[width=0.5\textwidth]{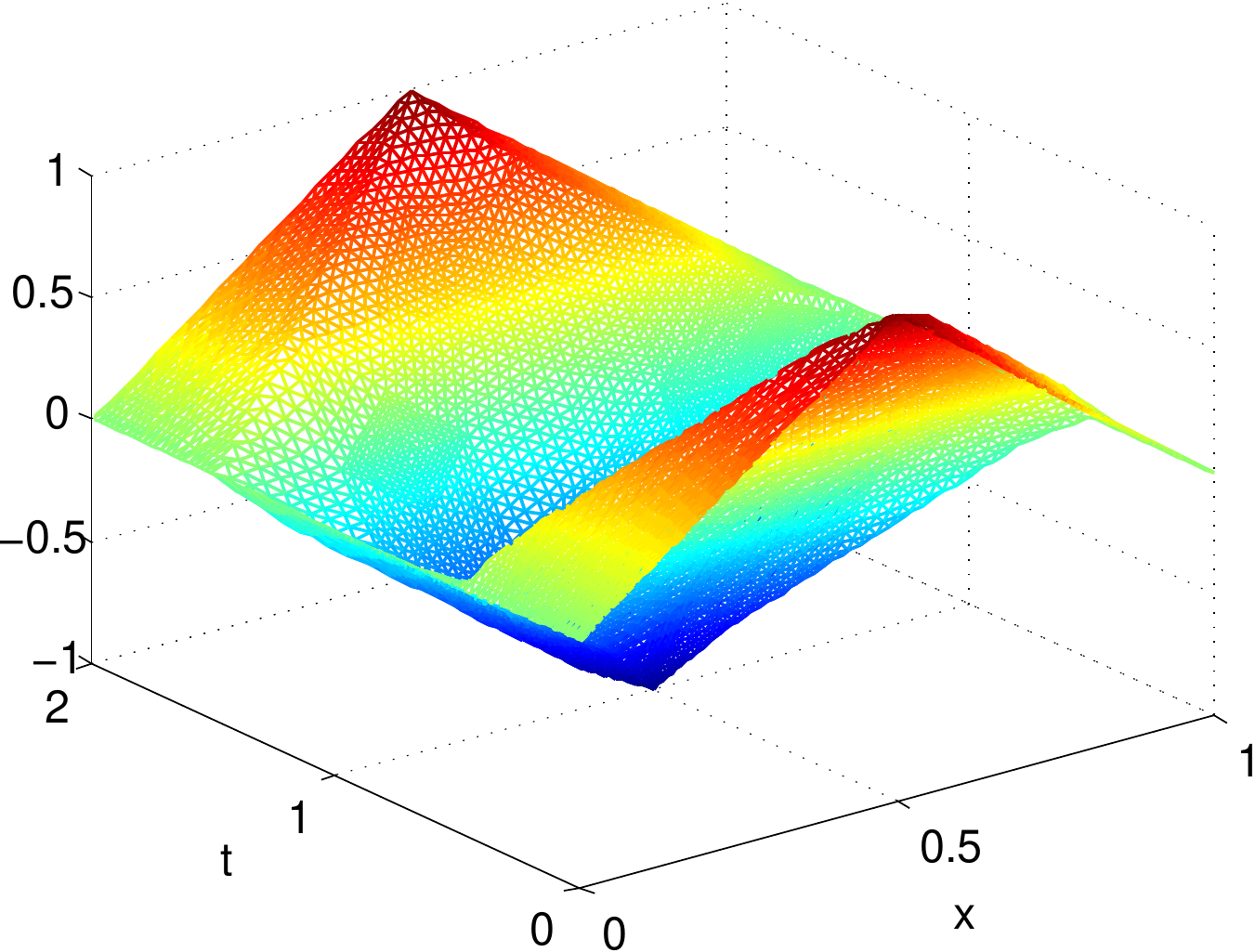} &
\includegraphics[width=0.5\textwidth]{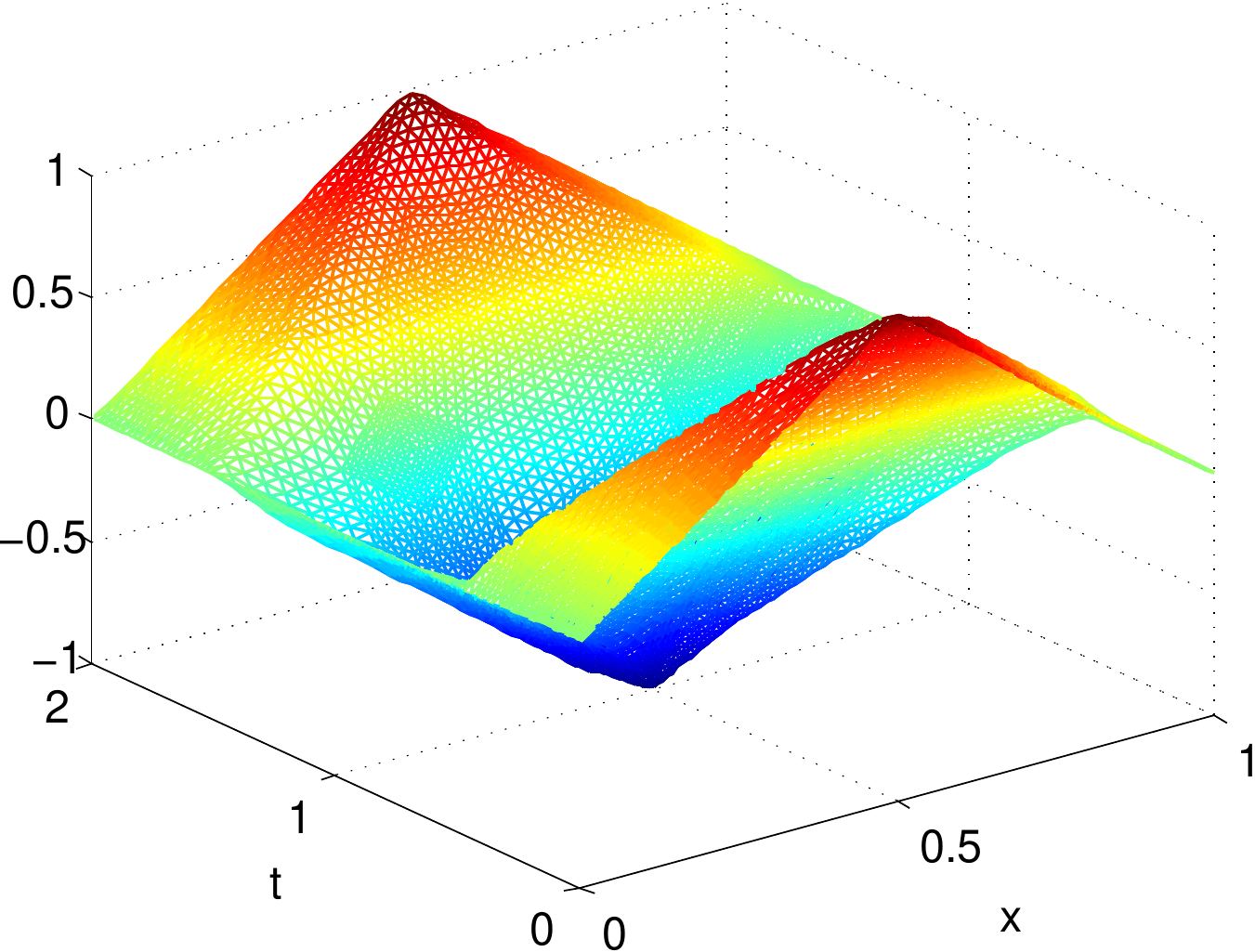} \\
(a) & (b)
\end{tabular}
\caption{Example \textbf{(EX2)} (a) Reference solution. (b) Solution reconstructed from the observation $y_{obs} = y|_{q_T^2}$.}
\label{fig:ex_vs_R}
\end{figure}

The exact solution $y$ corresponding to initial data \textbf{(EX2)} is displayed in Figure \ref{fig:ex_vs_R} (a) using the third mesh of the domain in Figure \ref{fig:qTvar} (b). 
Figure \ref{fig:ex_vs_R} (b) illustrates the solution $y_h$ of the mixed formulation \eqref{eq:mfh}, where the observation $y_{obs}$ is obtained as the restriction of $y$ to $q_T^2$.

\subsection{Two-dimensional space case ($N = 2$)}

We now illustrate the method introduced in Section \ref{recovering_y} in the two-dimensional case.  The procedure is similar but a bit more involved on a computational point of view, since $Q_T$ is now a subset of $\mathbb{R}^3$. 

\begin{figure}[ht]
\centering
\begin{tabular}{ccc}
\includegraphics[width=0.32\textwidth]{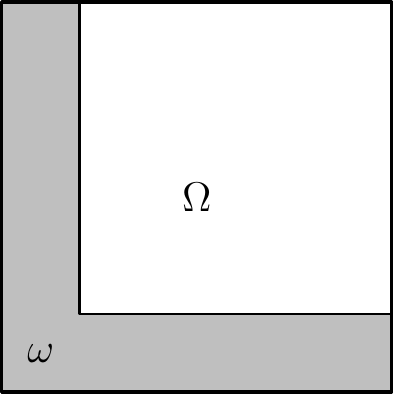} &\qquad&
\includegraphics[width=0.32\textwidth]{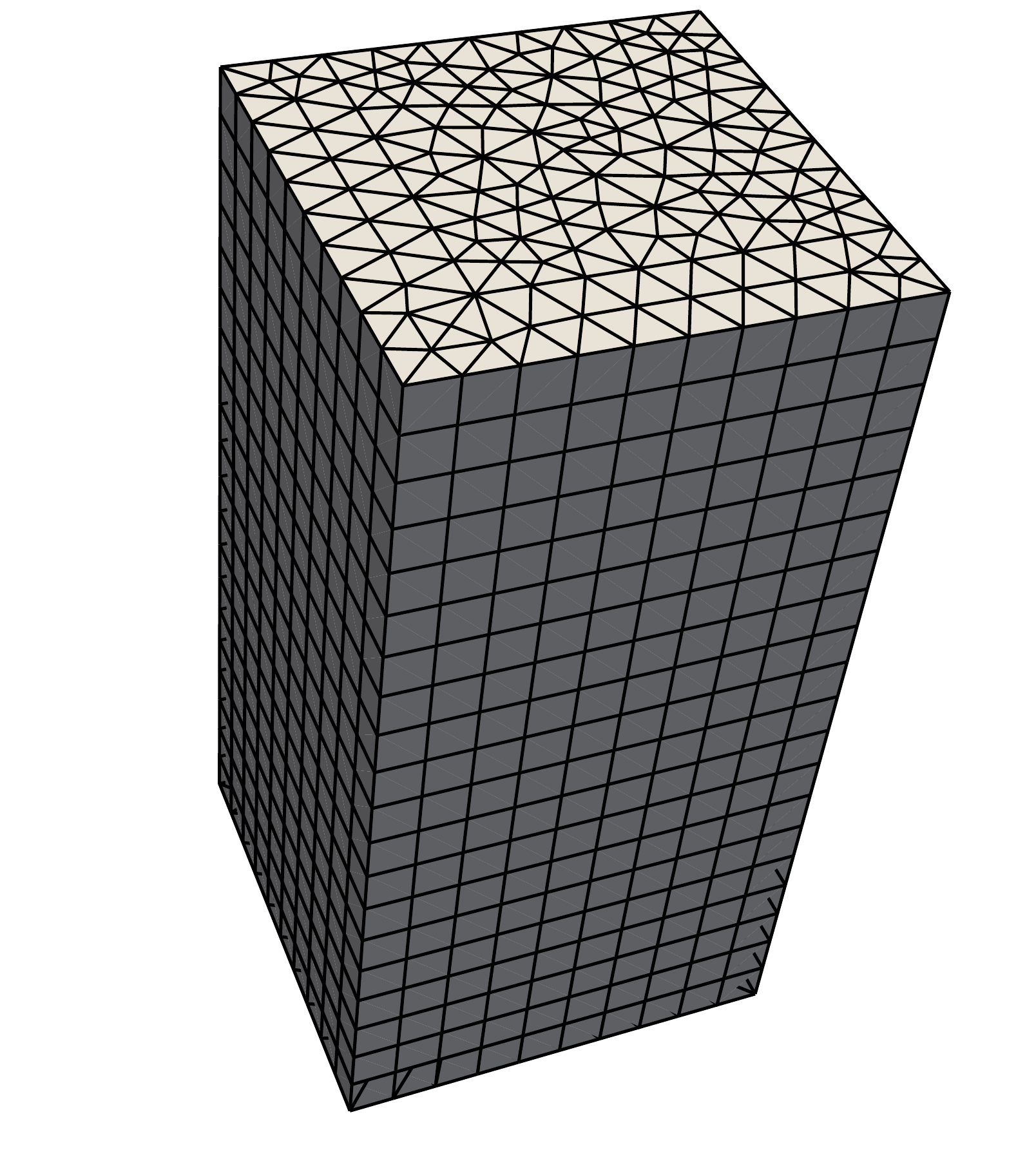} \\
(a) && (b)
\end{tabular}
\caption{(a) Example of sets $\Omega$ and $\omega$. (b) Example of mesh for $\Omega=(0,1)^2$ and $T = 2$.}
\label{fig:mesh}
\end{figure}

In order to approach the mixed-formulation \eqref{eq:mf}, we consider a mesh $\mathcal{T}_h$ of the domain $Q_T = \Omega \times (0,T)$ formed by triangular prisms. This mesh is obtained by extrapolating along the time axis a triangulation of the spatial domain $\Omega$. For an example in the case $\Omega=(0,1)^2$ and $T=2$ see Figure \ref{fig:mesh} (b) and for an example in the case of non-rectangular domains $\Omega \subset \mathbb{R}^2$ see Figure \ref{fig:heart} (b). For both examples, the extrapolation along the the time axis is uniform :  the height of the prismatic elements $\Delta t$ is constant.

Let $Z_h$ be the finite dimensional space defined as follows
\begin{equation}
\label{eq:Yh}
Z_h = \left\{ 
\begin{array}{l|l}
\ph_h = \psi(x_1, x_2) \theta(t) \in C^1(Q_T) & \psi|_{K_{x_1x_2}} \in \mathbb{P}(K_{xy}), \theta|_{K_t} \in \mathbb{Q}(K_t) \\ 
 \ph_h = 0 \text{ on } \Sigma_T & \text {for every } K = K_{x_1x_2} \times K_t \in \mathcal{T}_h.
\end{array}
\right\},
\end{equation}
$\mathbb{P}(K_{x_1x_2})$ is the space of functions corresponding to the reduced {\em Hsieh-Clough-Tocher} (HCT for short) $C^1$-element recalled in Section 4.1.1; $\mathbb{Q}(K_t)$ is a space of degree three polynomials on the interval $K_t$ of the form $[t_j,t_{j+1}]$ defined uniquely by their value and the value of their first derivative at the point $t_j$ and $t_{j+1}$. In other words, $Y_h$ is the finite element space obtained as a tensorial product between the reduced HCT finite element and cubic Hermite finite element. 
We check that on each element $K=K_{x_1x_2}\times K_t$, the function $\ph_h$ is
determined uniquely in term of the values of $\Sigma_K:=\{\ph(a_i),\ph_{x_1}(a_i),\ph_{x_2}(a_i),\ph_t(a_i),\ph_{x_1,t}(a_i),\ph_{x_2,t}(a_i), i=1,\cdots,6\}$ at the six nodes $a_i$ of $K$. Therefore, $\dim \Sigma_K=36$.

Similarly, let $\Lambda_h$ be the finite dimensional space defined by
\begin{equation}
\label{eq:Lambdah}
\Lambda_h = \left\{ 
\begin{array}{l|l}
\hspace{-0.1cm}\ph_h = \psi(x_1, x_2) \theta(t) \in C^0(Q_T) & \psi|_{K_{x_1x_2}} \in \mathbb{P}_1(K_{x_1x_2}), \theta|_{K_t} \in \mathbb{Q}_1(K_t) \\ 
 \ph_h = 0 \text{ on } \Sigma_T & \text {for every } K = K_{x_1x_2} \times K_t \in \mathcal{T}_h
\end{array}
\right\},
\end{equation}
where $\mathbb{P}_1(K_{x_1x_2})$ and  $\mathbb{Q}_1(K_t)$ are the spaces of degree one polynomials on the triangle $K_{x_1x_2}$ and interval $K_t$ respectively.

For any $h$, we check that $Z_h\subset Z$ and that $\Lambda_h\subset \Lambda$.

\subsubsection{Wave equation in a square}
\label{sec:square}
We first consider the case $\Omega$ defined by the unit square and again some explicit solutions used in \cite{cindea_moireau}.
Precisely, we define the following smooth initial condition:
\begin{equation}
(\textbf{EX1--2D}) \quad
\left\{
\begin{array}{l}
y_0(x_1, x_2) = 256 x_1^2 x_2^2 (1-x_1)^2 (1-x_2)^2 \\
y_1(x_1, x_2) = (1 - |2x_1 - 1|)(1 - |2x_2 - 1|)
\end{array}
\right. \quad (x_1, x_2) \in \Omega
\end{equation}
The corresponding solution of (\ref{eq:wave}) with $c\equiv 1, \ d\equiv 0$ and $f\equiv 0$ is given by :
\begin{equation}\label{eq:solex}
y(x_1, x_2, t) = \sum_{k,l > 0} \left( a_{kl} \cos (\mu_{kl} t) + \frac{b_{kl}}{\mu_{kl}} \sin(\mu_{kl} t) \right) \sin(k \pi x)\sin(l \pi y),
\end{equation}
where $\mu_{kl} = \pi \sqrt{k^2 + l^2}$ for every $k, l \in \mathbb{Z}^*$ and 
\begin{align*}
&a_{kl} = 2^{10} \frac{(\pi^2 k^2 - 12)(\pi^2 l^2 - 12)}{\pi^{10} k^5 l^5}((-1)^k - 1)((-1)^l - 1) \\
&b_{kl} = \frac{2^5}{\pi^4 k^2 l^2} \sin \frac{\pi k}{2} \sin \frac{\pi l}{2}.
\end{align*}
We also define the following initial data $(y_0, y_1) \in H_0^1(\Omega) \times L^2(\Omega)$:
\begin{equation}
(\textbf{EX2--2D}) \quad
\left\{
\begin{array}{l}
y_0(x_1, x_2) = (1 - |2x_1 - 1|)(1 - |2x_2 - 1|) \\
y_1(x_1, x_2) = \boldsymbol{1}_{(\frac{1}{3}, \frac{2}{3})^2}(x_1, x_2)
\end{array}
\quad (x_1, x_2) \in \Omega.
\right.
\end{equation} 
The Fourier coefficients of the corresponding solution are 
\begin{align*}
&a_{kl} = \frac{2^5}{\pi^4 k^2 l^2} \sin \frac{\pi k}{2} \sin\frac{\pi l}{2} \\
&b_{kl} = \frac{1}{\pi^2 k l} \left( \cos \frac{\pi k}{3} - \cos \frac{2\pi k}{3}\right) \left( \cos \frac{\pi l}{3} - \cos \frac{2\pi l}{3}\right).
\end{align*}

In what follows, we consider $\omega$ the subset of $\Omega$ described in Figure \ref{fig:mesh} (a) and given by:
\begin{equation}
\label{eq:omegaS}
\omega = \left( (0, 0.2) \times (0,1) \right) \cup \left( (0,1) \times (0, 0.2) \right).
\end{equation}

It is easy to see that this choice of $\omega$ and $T = 2$ provide a domain $q_T=\omega\times (0,T)$ which satisfies the geometric optic condition, and, hence, inequality \eqref{iobs} holds. We consider 3 levels of meshes of $Q_T$, labelled from 1 to 3 and containing the number of elements (prisms) and nodes listed in Table \ref{tab:mesh}.

\begin{table}[ht]
\centering
\begin{tabular}{|c|ccc|}
\hline 
Mesh Number  & 1 & 2 & 3 \\
\hline
Number of elements & 5 320 & 15 320 & 42 230 \\
Number of nodes & 3 234 & 8 799 & 23 370\\
$\Delta t$ & $0.2$ & $0.1$ & $0.05$ \\
\hline
\end{tabular}
\caption{Characteristics of the meshes used for $Q_T = (0,1)^2 \times (0,2)$.}
\label{tab:mesh}
\end{table}

For each of these meshes we solve the mixed formulation \eqref{eq:mf} with the term $y_{obs}$ appearing in the right-hand side obtained as the restriction to $q_T$ of the solution computed by \eqref{eq:solex} for initial data \textbf{EX1--2D} and \textbf{EX2--2D}. 

Table \ref{tab:EX1} concerns the example \textbf{EX1--2D}. In this table we list the norm of the relative error between the exact solution $y$ given by \eqref{eq:solex} and the solution $y_h$ of the mixed formulation \eqref{eq:mf}, the $L^2$ norm of $Ly_h$ and the $L^2$ norm of the Lagrange multiplier $\lambda_h$.
\begin{table}[ht]
\centering
\begin{tabular}{|c|ccc|}
\hline
Mesh number & 1 & 2 & 3 \\
\hline
$\frac{\|y-y_h\|_{L^2(Q_T)}}{\|y\|_{L^2(Q_T)}}$ & $4.58 \times 10^{-2}$ & $3.18 \times 10^{-2}$ & $1.38 \times 10^{-2}$ \\
$\|Ly_h\|_{L^2(Q_T)}$ & $1.44$ & $1.05$ & $1.05$\\
$\|\lambda_h\|_{L^2(Q_T)}$ & $2.87 \times 10^{-5}$ & $1.36 \times 10^{-5}$ & $7.34 \times 10^{-6}$ \\
$\sharp$ CG iterates & 121 & 180 & 168\\
\hline
\end{tabular}
\caption{$\varepsilon = 0$: Example \textbf{EX1--2D} - $r = 1$.}
\label{tab:EX1}
\end{table}

As theoretically stated in Remark \ref{rk_lambda_sys} and observed in numerical experiments in the case $N = 1$ (see, for instance, Table \ref{tab:EX1}), the Lagrange multiplier $\lambda_h$ vanishes as $h\to 0$. In Table \ref{tab:EX1-lambdazero} we display the results obtained by numerically solving the variational problem \eqref{eq:mf} obtained from the mixed formulation when $\lambda_h = 0$.

\begin{table}[ht]
\centering
\begin{tabular}{|c|ccc|}
\hline
Mesh number & 1 & 2 & 3 \\
\hline
$\frac{\|y-y_h\|_{L^2(Q_T)}}{\|y\|_{L^2(Q_T)}}$ & $7.05 \times 10^{-2}$ & $4.44 \times 10^{-2}$ & $2.37 \times 10^{-2}$ \\
$\|Ly_h\|_{L^2(Q_T)}$ & $1.31$ & $0.97$ & $0.97$ \\
\hline
\end{tabular}
\caption{Example \textbf{EX1--2D} -- $r = 1$ -- $\lambda$ fixed to zero.}
\label{tab:EX1-lambdazero}
\end{table}

Tables \ref{tab:EX2} and \ref{tab:EX2-lambdazero} display the results obtained for the initial data specified by \textbf{EX2--2D}, for the solutions $(y_h, \lambda_h)$ of the mixed formulation and for the variational problem obtained when $\lambda_h = 0$ respectively.

\begin{table}[ht]
\centering
\begin{tabular}{|c|ccc|}
\hline
Mesh number & 1 & 2 & 3 \\
\hline
$\frac{\|y-y_h\|_{L^2(Q_T)}}{\|y\|_{L^2(Q_T)}}$ & $4.74 \times 10^{-2}$ & $3.72 \times 10^{-2}$ & $2.09 \times 10^{-2}$  \\
$\|Ly_h\|_{L^2(Q_T)}$ & $1.18$ & $0.89$ & $1.06$\\
$\|\lambda_h\|_{L^2(Q_T)}$ & $3.21 \times 10^{-5}$ & $1.46 \times 10^{-5}$ & $1.17 \times 10^{-5}$ \\
$\sharp$ CG iterates & $128$ & $191$ & $168$ \\
\hline
\end{tabular}
\caption{Example \textbf{EX2--2D} -- $r = 1$.}
\label{tab:EX2}
\end{table}

The results are similar for both examples. In both cases we observe a linear convergence of $y_h$ to $y$ in the norm $L^2$ over $Q_T$ when $h$ goes to zero. Similarly, the norm $\|\lambda_h\|_{L^2(Q_T)}$ linearly decreases as $h$ goes to zero. 

\begin{table}[ht!]
\centering
\begin{tabular}{|c|ccc|}
\hline
Mesh number & 1 & 2 & 3 \\
\hline
$\frac{\|y-y_h\|_{L^2(Q_T)}}{\|y\|_{L^2(Q_T)}}$ & $6.75 \times 10^{-2}$ & $4.93 \times 10^{-2}$ & $3.37 \times 10^{-2}$  \\
$\|Ly_h\|_{L^2(Q_T)}$ & $1.07$ & $0.82$ & $0.97$ \\
\hline
\end{tabular}
\caption{Example \textbf{EX2--2D} -- $r = 1$ -- $\lambda$ fixed to zero.}
\label{tab:EX2-lambdazero}
\end{table}


\subsubsection{Wave equation in a non-rectangular domain of $\mathbb{R}^2$}

\begin{figure}[ht]
\hspace{-0.6cm}
\centering
\begin{tabular}{cc}
\includegraphics[width=0.51\textwidth]{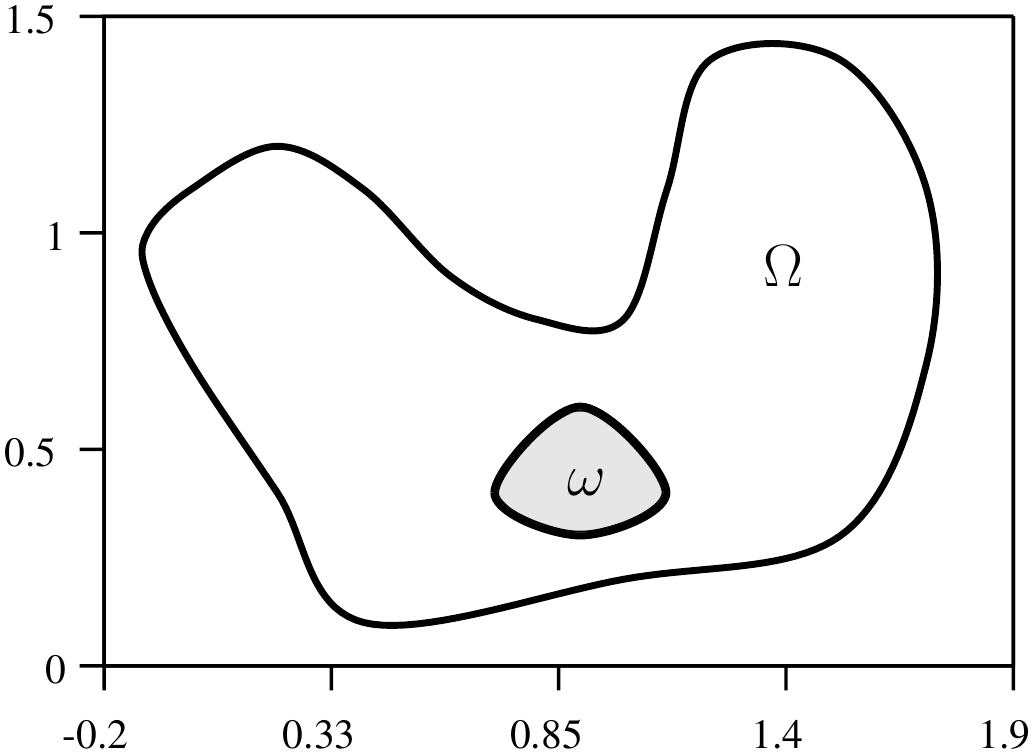} & 
\includegraphics[width=0.44\textwidth]{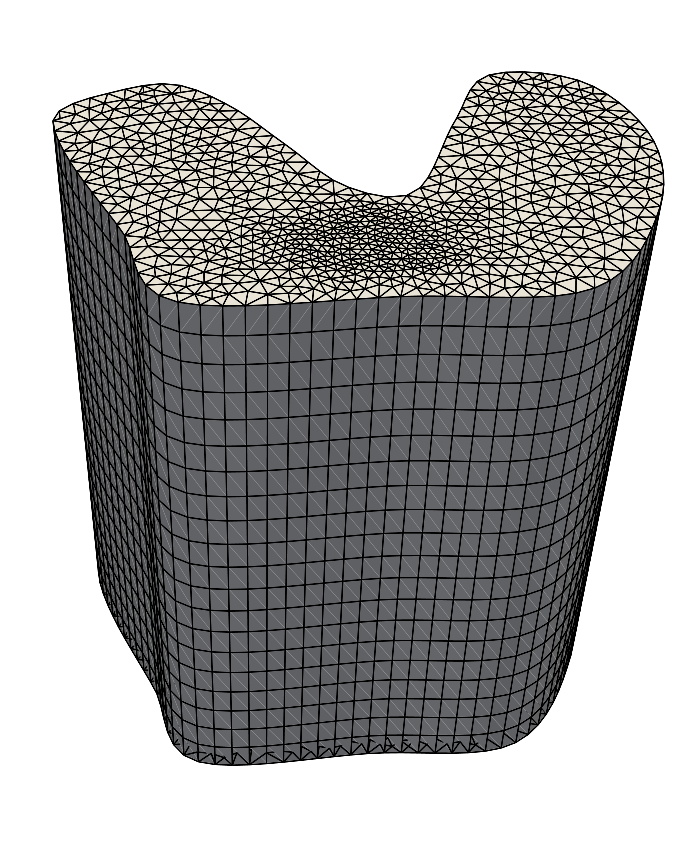} \\
(a) & (b)
\end{tabular}
\caption{(a) Example of sets $\Omega$ and $\omega$. (b) Example of mesh of the domain $Q_T$.}
\label{fig:heart}
\end{figure}

Let $\Omega \subset \mathbb{R}^2$ be a domain with a regular boundary and $\omega$ a non-empty subset with regular boundary. An example of such a configuration is illustrated in Figure \ref{fig:heart} (a). As in the previous section, we take  $T = 2$ and we build a mesh formed by triangular prisms of the domain $Q_T = \Omega \times (0, T)$. An example of such a mesh associated to the domain $\Omega$ is displayed in Figure \ref{fig:heart} (b). This mesh is composed by $17 \ 934$ nodes distributed in $32 \ 140$ prismatic elements (this mesh corresponds to the mesh number 2 described in Table \ref{tab:mesh_h}).

We consider three levels of meshes of the domain $Q_T$ formed by the number of prisms and containing the number of nodes reported in Table \ref{tab:mesh_h}.

\begin{table}[ht]
\centering
\begin{tabular}{|c|ccc|}
\hline
Mesh number & 1 & 2 & 3\\
\hline 
Number of elements & 5  730 & 32 1400 & 130 280 \\
Number of nodes & 3 432 & 17 934 & 69 864 \\
Height of elements ($\Delta t$) & 0.2 & 0.1 & 0.05 \\
\hline
\end{tabular}
\caption{Characteristics of the three meshes associated with $Q_T$.}
\label{tab:mesh_h}
\end{table}

Comparing to the situation described in Subsection \ref{sec:square}, the eigenfunctions and eigenvectors of the Dirichlet Laplace operator defined on $\Omega$ are not explicitly available here. Consequently, from a given set of initial data, we numerically solve the wave equation \eqref{eq:wave} using a standard time-marching method, from which we can extract an observation on $q_T$. Precisely, we use a $P_1$ finite elements method in space coupled with a Newmark unconditionally stable scheme for the time discretization. Hence, we solve the wave equation on  the same mesh which was extrapolated in time in order to obtain the mesh number 2 of $Q_T$. This two-dimensional mesh contains $1 \ 704$ nodes and $3 \ 257$ triangles. The time discretization step is $\Delta t = 10^{-2}$. We denote $\overline{y}_h$ the solution obtained in this way for the initial data $(y_0, y_1) \in H_0^1(\Omega) \times L^2(\Omega)$ given by

\begin{equation}
\label{eq:idheart}
\left\{ 
\begin{array}{ll}
-\Delta y_0 = 10,  &\quad \text{in } \Omega \\
y_0 = 0, & \quad \text{on } \partial \Omega,
\end{array}
\right.
\qquad y_1 =0.
\end{equation}

From $\overline{y}_h$ we generate the observation $y_{obs}$ as the restriction of $\overline{y}_h$ to $q_T$. Finally, from this observation we reconstruct $y_h$ as the solution of the mixed formulation \eqref{eq:mfeps} on each of the three meshes described in Table \ref{tab:mesh_h}.  Table \ref{tab:h} display some norms of $y_h$ and $\lambda_h$ obtained for the three meshes and illustrates again the convergence of the method.

\begin{table}[ht!]
\centering
\begin{tabular}{|c|ccc|}
\hline
Mesh number & 1 & 2 & 3\\
\hline
$\frac{\|\overline{y}_h - y_h\|_{L^2(Q_T)}}{\|\overline{y}_h\|_{L^2(Q_T)}}$ & $1.88 \times 10^{-1}$ & $8.04 \times 10^{-2}$ & $7.11 \times 10^{-2}$ \\
$\|Ly_h\|_{L^2(Q_T)}$ & $3.21$ & $2.01$ & $1.57$\\
$\|\lambda_h\|_{L^2(Q_T)}$ & $8.26 \times 10^{-5}$ & $3.62 \times 10^{-5}$ & $2.84 \times 10^{-5}$ \\
$\sharp$ CG iterates & 52 & 167 &  400 \\
\hline
\end{tabular}
\caption{Initial data $(y_0,y_1)$ given by (\ref{eq:idheart}) -  $r = 1$.}
\label{tab:h}
\end{table}

\begin{figure}[ht]
\begin{tabular}{cc}
\includegraphics[width=0.5\textwidth]{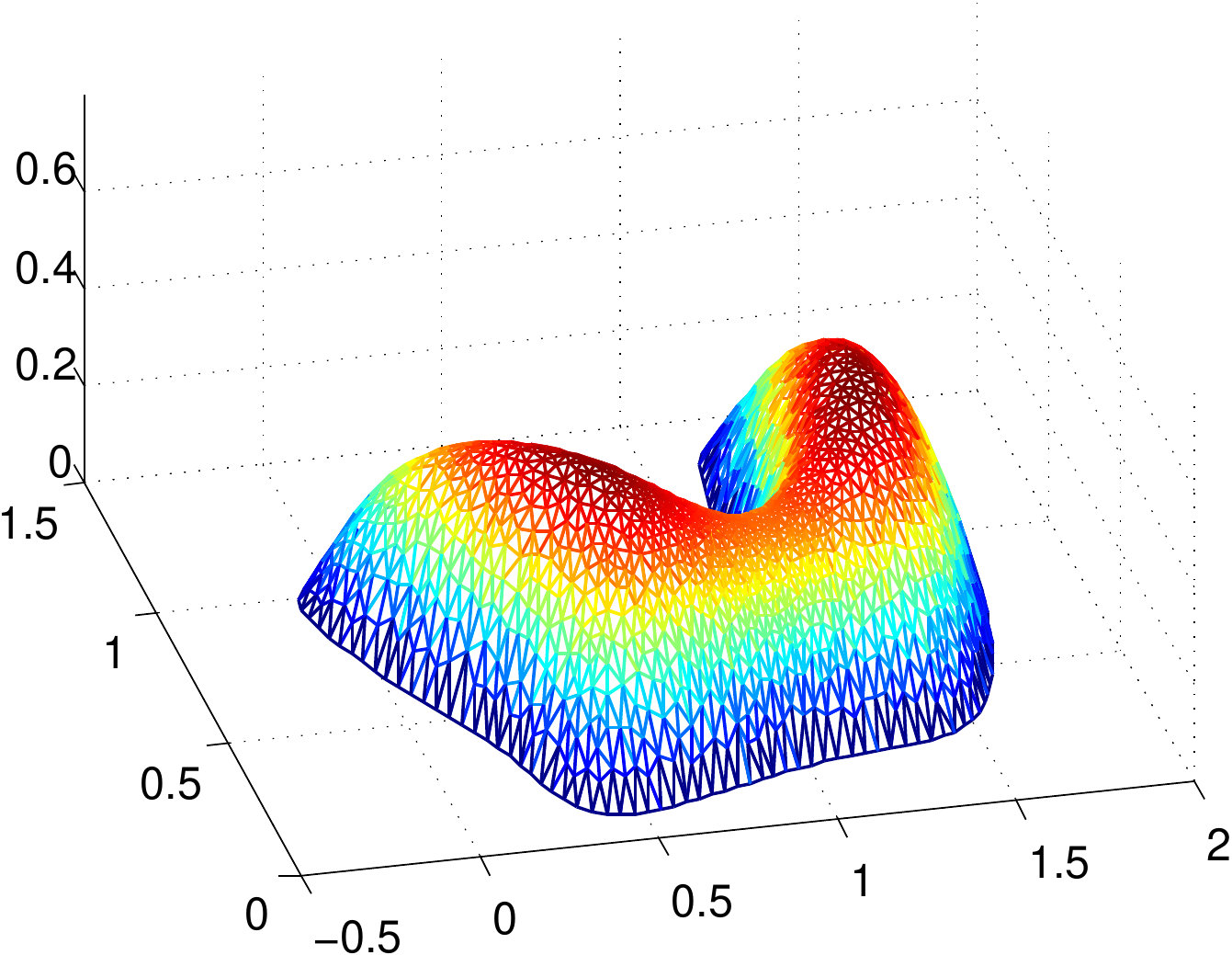} & \includegraphics[width=0.5\textwidth]{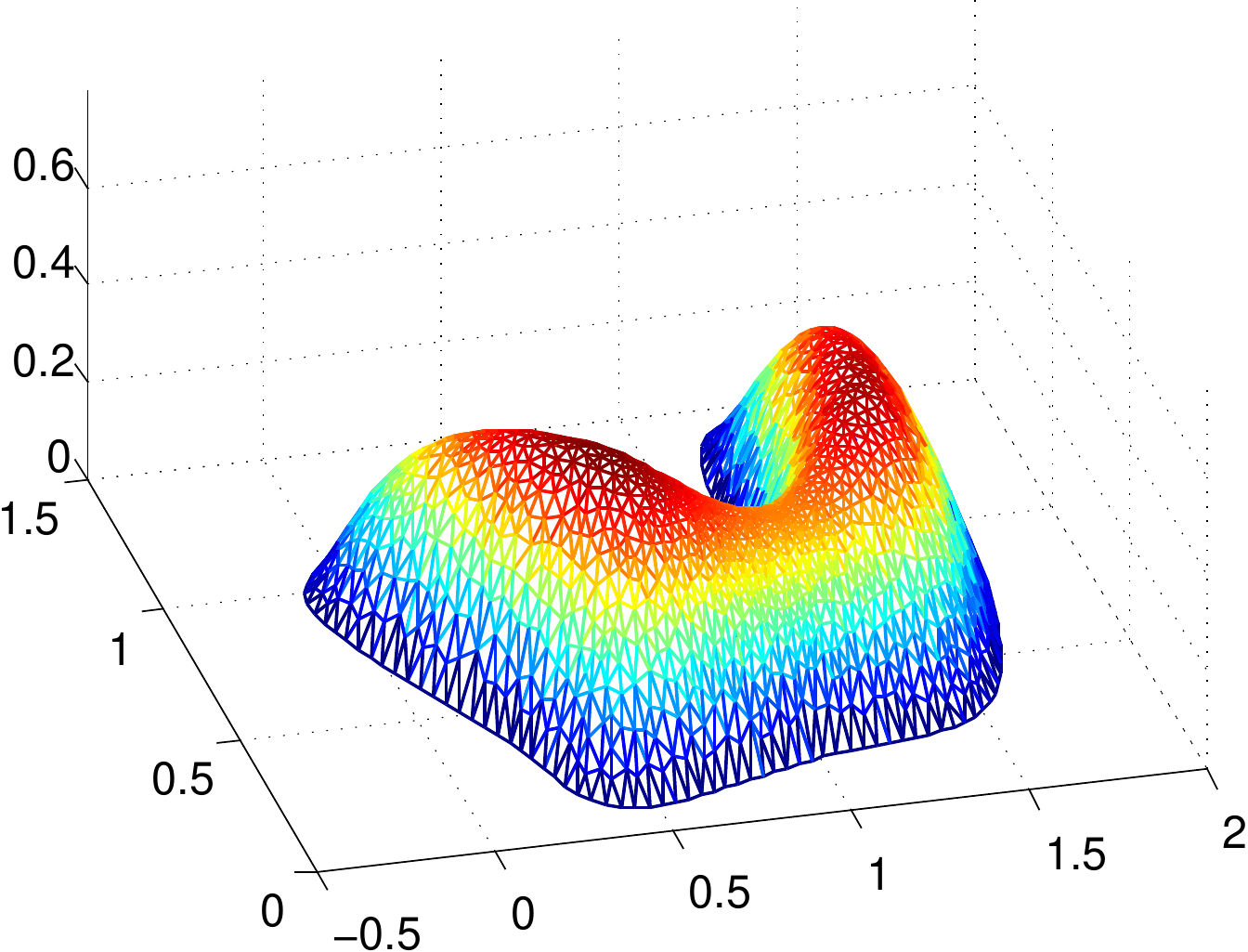}  \\
(a) & (b)
\end{tabular}
\caption{(a) Initial data $y_0$ given by \eqref{eq:idheart}. (b) Reconstructed initial data $y_h(\cdot, 0)$.}
\label{fig:idheart}
\end{figure}

Figure \ref{fig:idheart} (a) displays the solution $y_0$ of \eqref{eq:idheart} and Figure \ref{fig:idheart} (b) displays the initial position $y_h(\cdot, 0)$ corresponding to the solution of our inverse problem. The error between these two functions is given by $\|y_0 - y_h(\cdot, 0)\|_{L^2(\Omega)} = 2.05 \times 10^{-2}$ which is consistent with the results reported in Table \ref{tab:h}.

\section{Concluding remarks and perspectives}  \label{sec_conclusion}

The mixed formulations we have introduced here in order to address inverse problems for the wave equation seems original. These formulations are nothing else than the Euler systems associated to least-squares type functionals and depend on both the state to be reconstruct and a Lagrange multiplier. This Lagrange multiplier is introduced to take into account the state constraint $Ly-f=0$ and turns out to be the controlled solution of a wave equation with the source term $(y-y_{obs})\,1_{q_T}$. This approach, recently used in a controllability context in \cite{NC-AM-mixedwave}, leads to a variational problem defined over time-space functional Hilbert spaces, without distinction between the time and the space variable. 
The main ingredient allowing to prove the well-posedness of the mixed formulation and therefore the reconstruction of the solution, is a generalized observability inequality, assuming here 
some geometric conditions on the observation zone. 

At the practical level, the discrete mixed time-space formulation is solved in a systematic way in the framework of the finite element theory. The approximation is conformal allowing to obtain the strong convergence of the approximation as the discretization parameters tends to zero. In particular, we emphasize that there is no need, contrarily to the classical approach, to prove some uniform discrete observability inequality: we simply use the observability equality on the finite dimensional discrete space. The resolution amounts to solve a sparse symmetric linear system : the corresponding matrix can be preconditioned if necessary, and may be computed once for all as it does not depend on the observation $y_{obs}$. Eventually, the space-time discretization of the domain allows an adaptation of the mesh so as to reduce the computational cost and capture the main features of the solutions.  Similarly, this space-time formulation is very appropriate to the non-cylindrical situation.

In agreement with the theoretical convergence, the numerical experiments reported here display a very good behavior and robustness of the approach: the reconstructed approximate solution converges strongly to the solution of the wave equation associated to the available observation. Remark that from the continuous dependence of the solution with respect to the observation, the method is robust with respect to the possible noise on the data. 

As mentioned at the end of Section \ref{recovering_y_f}, additional assumption on the source term allows to determine uniquely the pair $(y,f)$ from a partial measurement on $q_T$ or on a part $\Sigma_T$ sufficiently large of the boundary. For instance, from \cite[Theorem 2.1]{yamamoto}, assuming that the source term takes the form $f(x,t)=\sigma(t) \mu(x)$ with $\sigma\in C^1([0,T])$, $\sigma(0)\neq 0$ and $\mu\in H^{-1}(\Omega)$, then the following holds: there exists a positive constant $C$ such that 
\begin{equation}
\Vert \mu\Vert^2_{H^{-1}(\Omega)} \leq C\biggl( \biggl\Vert \frac{\partial y}{\partial \nu}\biggr\Vert^2_{L^2(\Sigma_T)} + \Vert Ly-\sigma(t)\mu(x)\Vert^2_{L^2(Q_T)}\biggr), \quad \forall (y,\mu)\in S
\end{equation}
where $y$ solves (\ref{eq:wave}) with $(y_0,y_1)\equiv 0$, $c=1$ and $(\Sigma_T,T,Q_T)$ satisfies a geometric condition and $S$ denotes an appropriate functional space.  Using this inequality (similar to \ref{iobs}), we can study the mixed formulation associated to the Lagrangian from $S\times L^2(Q_T)\to \mathbb{R}$ defined by
$$
\mathcal{L}((y,\mu),\lambda):= \frac{1}{2}\biggl\Vert \frac{\partial y}{\partial\nu}-y_{obs}\biggr\Vert^2_{L^2(\Sigma_T)} + \int_{Q_T} \lambda (Ly-\sigma \mu) \, dx\,dt
$$
to fully reconstruct $y$ and $\mu$ from $y_{obs}$ and $\sigma$.

Eventually, since the mixed formulations rely essentially on a generalized observability inequality, it may be employed to any other observable systems for which such property is available :
we mention notably the parabolic case usually -- in view of regularization property -- badly conditioned and for which direct and robust methods are certainly very advantageous. We refer to \cite{munch_desouza}  where this issue is investigated.


\bibliographystyle{siam}




\end{document}